\documentclass[11pt]{article}
\newtheorem{theorem}{Theorem}[section]

\newtheorem{proposition}[theorem]{Proposition}
\newtheorem{corollary}[theorem]{Corollary}
\newtheorem{lemma}[theorem]{Lemma}
\newtheorem{example}[theorem]{Example}

\def\qed{\hfill $\Box$ \medskip}
\usepackage{amsmath}
\usepackage{amssymb}
\usepackage{graphicx}
\def\IR{{\bf R}}
\def\IC{{\bf C}}
\def\IF{{\bf F}}
\def\cS{{\cal S}}
\def\cR{{\cal R}}
\def\cC{{\cal C}}
\def\cF{{\cal F}}
\def\bK{{\bf K}}
\def\bC{{\bf C}}
\def\a{{\alpha}}
\def\b{{\beta}}
\def\t{{\theta}}

\def\tc{{\theta_C}}
\def\td{{\theta_D}}
\def\tp{{\theta_P}}
\def\arg{{\rm arg}\,}
\def\Ra{{\ \Rightarrow\ }}
\def\Lra{{\ \Leftrightarrow\ }}
\def\logm{{\prec_{\rm log}}\,}
\def\diag{{\rm diag}\,}
\def\tr{{\rm tr}\,}
\def\rank{{\rm rank}\,}
\def\span{{\rm span}\,}
\def\ker{{\rm Ker}\,}
\def\[{\left [}
\def\]{\right ]}
\def\({\left (}
\def\){\right )}
\def\arg{{\rm arg}}
\def\dfrac{\displaystyle\frac}

\usepackage{tikz}
\usetikzlibrary{patterns}

\begin{document}
\openup .38\jot
\title{Minkowski product of convex sets and product numerical range}

\author{Chi-Kwong Li\thanks{Department of Mathematics, College of William and Mary,
Williamsburg,
VA 23187, USA. ckli@math.wm.edu}, Diane Christine Pelejo\thanks{Department of Mathematics,
College of William and Mary, Williamsburg, VA 23187, USA. dppelejo@wm.edu}, Yiu-Tung
Poon\thanks{ Department of Mathematics, Iowa State University, Ames, IA 50011,
USA. ytpoon@iastate.edu},
Kuo-Zhong Wang\thanks{Department of Applied Mathematics,
National Chiao Tung University, Hsinchu 30010, Taiwan.
kzwang@math.nctu.edu.tw}}
\date{}
\maketitle

\centerline{In memory of Leiba Rodman.}

\begin{abstract}
Let $K_1, K_2$ be two compact convex sets in $\IC$. Their Minkowski product is the set
$K_1K_2 = \{ab: a \in K_1, b\in K_2\}$. We show that the set $K_1K_2$ is star-shaped
if $K_1$ is a line segment or a circular disk.
Examples for $K_1$ and $K_2$ are given so that $K_1$ and $K_2$ are triangles (including
interior) and $K_1K_2$ is not star-shaped. This gives a negative answer to a conjecture
by  Puchala et. al concerning the product numerical range in the study of quantum information
science. Additional results and open problems are presented.
\end{abstract}



\section{Introduction}

Let $K_1, K_2$ be compact convex sets in $\IC$. We study
the Minkowski product of the sets defined and denoted by
$$K_1K_2 = \{ab:a \in K_1, b \in K_2\}.$$
This topic arises naturally in many branches of research. For example,
in numerical analysis, computations are subject to errors
caused by the precision of the machines and round-off errors.
Sometimes measurement errors in the raw data may also affect the
accuracy. So, when two real numbers $a$ and $b$ are multiplied, the actual
answer may actually be the product of numbers in two intervals containing $a$ and $b$;
when two complex numbers $a$ and $b$ are multiplied, the actual answer may actually be
the product of numbers from two regions in the complex plane. The study of the product set
also has applications in computer-aided design, reflection and refraction of wavefronts in
geometrical optics, stability
characterization of multi-parameter control systems, and the shape analysis and
procedural generation of two-dimensional domains. For more discussion about
these topics, see \cite{FMR} and the references therein.
Another application comes from the study of quantum information science.
For a complex $n\times n$ matrix $A$, its numerical range is defined and denoted by
$$W(A) = \{ x^*Ax: x \in \IC^n, x^*x = 1\}.$$
The numerical range of a matrix is always a compact convex
set and carries a lot of information about the matrix, e.g., see \cite{HJ}.

Denote by $X \otimes Y$ the Kronecker product of two matrices or vectors.
Then the decomposable numerical range of $T \in M_m \otimes M_n \equiv M_{mn}$ is defined by
$$W^{\otimes}(T) = \{(x^*\otimes y^*)T(x\otimes y):
x \in \IC^m, y \in \IC^n, x^*x = y^*y = 1\},$$
which is a subset of $W(T)$. In the context of quantum information science,
this set corresponds to the collection of $\langle T, P \otimes Q\rangle$,
where $P \in M_m, Q \in M_n$ are pure states (i.e., rank one orthogonal projections).
In particular, if
$T = A \otimes B$ with $(A,B) \in M_m \times M_n$, then
$$W^{\otimes}(A\otimes B) = \{(x^*\otimes y^*)(A\otimes B)(x\otimes y):
x \in \IC^m, y \in \IC^n, x^*x = y^*y = 1\} = W(A)W(B).$$
So, the set $W^{\otimes}(A\otimes B)$ is just the Minkowski product of the two compact
convex sets $W(A)$ and $W(B)$.
In particular, the following was proved in  \cite{Karol}. (Their proofs concern the
product numerical range that can be easily adapted to general compact convex sets.)

\begin{proposition} \label{1.0} Suppose $K_1, K_2$ are compact convex sets in $\IC$.

{\rm (a)}
The set
$K_1 K_2$ is simply connected.

{\rm (b)}
If $0 \in K_1 \cup K_2$, then $K_1K_2$ is star-shaped with $0$ as a star center.

\end{proposition}

It was conjectured in \cite{Karol} that
the set $K_1K_2$ is always star-shaped.
In this paper, we will show that the conjecture is not true in general (Section 3.1).
The proof depends on a detailed analysis of the product sets of two closed line
segments (Section 2). Then we obtain some
conditions under which the product set of two convex polygons is star-shaped
(Sections 3.2). Furthermore, we show that $K_1K_2$ is star-shaped for any compact convex
set $K_2$ if $K_1$ is a closed line segment or a closed circular disk in Sections
4 and 5. Some additional results and open problems are mentioned
in Section 6. In particular, in Theorem \ref{6.2}, we will improve the
following result, which is a consequence of the simply connectedness of $K_1K_2$
\cite[Proposition 1]{Karol}.

\begin{proposition} \label{1.1} Suppose $K_1$, $K_2$ are compact convex
sets in  $\IC$ and $p\in K_1K_2$. 
Then $K_1K_2$ is  star-shaped with $p$ 
as a star center if and only if $K_1K_2$
contains the line segment joining $p$ 
to $ab$ for any $a \in \partial K_1$
and $b \in \partial K_2$.
\end{proposition}

In our discussion, we
denote by $K(z_1,z_2, \dots, z_m)$ the convex hull of the set
$\{z_1, \dots, z_m\} \subseteq \IC$. In particular, $K(z_1,z_2)$ is the line segment in $\IC$
joining $z_1, z_2$. Also, if $K_1 = \{\alpha\}$, we write $K_1K_2 = \alpha K_2$.

\section{The product set of two segments}

We first give a complete description of the set $K_1K_2$ when
$K_1 = K(\a_1, \a_2)$ and $K_2 = K(\b_1,\b_2)$ are two line segments. McAllister  has plotted
some examples in \cite{Mc} but the analysis is not complete.
In the context of product numerical range, it is known, see for example,
\cite[Theorem 4.3]{Li}, that
$W(T)$ is a line segment if and only if $T$ is normal with collinear eigenvalues.
In such a case, $W(T) = W(T_0)$ for a normal matrix $T_0 \in M_2$ having the two endpoints of
$W(T)$ as its eigenvalues.
Thus, the study of $K_1K_2$ when $K_1, K_2$ are close line segments corresponds to the
study of $W^{\otimes}(A\otimes B) = W(A)W(B)$ for $A \in M_m, B\in M_n$ with special structure,
and $W^{\otimes}(A\otimes B) = W^{\otimes}(A_0\otimes B_0)$ for some
normal matrices $A_0, B_0 \in M_2$.
We have the following result.

\begin{theorem}\label{2.1}
Let $K_1 = K(\a_1, \a_2)$ and $K_2 = K(\b_1,\b_2)$ be two line segments in
$\IC$. Then $K_1K_2$ is a star-shaped subset of
$K(\a_1\b_1,\a_1\b_2,\a_2\b_1,\a_2\b_2)$.
\end{theorem}

In general, $K(\alpha_1,\ldots, \alpha_n)K(\beta_1,\ldots, \beta_m) \subseteq K(\alpha_1\beta_1, \alpha_1,\beta_2,\ldots, \alpha_i\beta_j, \ldots, \alpha_n\beta_m)$
because
$$ \left(\sum\limits_{i}p_i\alpha_i\right)\left(\sum\limits_{j}q_j\beta_j\right)
=\left(\sum\limits_{i,j}p_iq_j\alpha_i\beta_j\right) $$
and $\sum_i p_i=1$ and $\sum_j q_j=1$ imply that $\sum\limits_{i,j}p_iq_j=1$. The key point of Theorem \ref{2.1} is the star-shapedness of the product of two line segments in $\IC$.

We will give a complete description of the set $K_1K_2$ in the following.
If one or both of the line segments $K_1, K_2$ lie(s) in a line passing through origin,
the description is relatively easy as shown in the following.

\medskip\noindent
\begin{proposition}\label{2.2} Let $K_1 = K(\a_1, \a_2)$ and $K_2 = K(\b_1,\b_2)$ be two line
segments in $\IC$.

\begin{itemize}
\item[{\rm 1.}] If both  $K(0, \a_1, \a_2)$  and $K(0, \b_1, \b_2)$ are line segments, then
$K_1 K_2$ is the line segment
$$K(\a_1\b_1,\a_1\b_2,\a_2\b_1,\a_2\b_2).$$

\item
[{\rm 2.}] Suppose $K(0, \a_1, \a_2)$  is a line segment and $K(0, \b_1, \b_2)$ is not.

{\rm (2.a)} If
$0 \in K(\a_1,\a_2)$, then $K_1K_2 = K(0, \a_1\b_1, \a_1\b_2) \cup K(0,\a_2\b_1, \a_2\b_2)$
is the union of two

\  triangles (one of them may degenerate to $\{0\}$)
meeting at 0, which is the star center of $K_1K_2$.

{\rm (2.b)} If $0 \notin K(\a_1,\a_2)$ then $K_1K_2 =  K(\a_1\b_1,\a_1\b_2,\a_2\b_1,\a_2\b_2)$.
\end{itemize}
\end{proposition}

\it Proof. \rm
{\rm 1.} There exist $\a,\ \b, a_1, a_2, b_1, b_2\in\IR$ such that
$K_1=\{re^{i\alpha}: a_1 \leq r\leq b_1\}$ and
$K_2=\{re^{i\beta}: a_2 \leq r\leq b_2\}$. So, we have
$$K_1K_2=\{re^{i(\alpha+\beta)}: a_3 \leq r\leq b_3\} \ \mbox{ for some }a_3,\ b_3\in \IR.$$

{\rm (2.a)} Evidently, $K_1K_2 =K(0,\a_1)K_2\cup K(0,\a_2)K_2$ and $K(0,\a_i)K_2\subseteq
K(0, \a_i\b_1, \a_i\b_2)$ for $i=1,\ 2$. We are going to show that $K(0,\a_i)K(\b_1,\b_2)
=K(0, \a_i\b_1, \a_i\b_2)$ for $i=1,\ 2$.

Clearly, $0\in K(0,\alpha_i)K(\beta_1,\beta_2)$. If $x\in K(0,\alpha_i\beta_1,\alpha_i\beta_2)\setminus\{0\}$, then there exist $s,\ t \ge 0$ with $0<s+t\le 1$ such that
$x=s \alpha_i\beta_1+t\alpha_i\beta_2$. Therefore,
 $x=ab$, where
$$a=(s+t)\alpha_i \in  K(0,\alpha_i) \quad \mbox{ and }\quad b=\frac{s}{s+t}\beta_1+\frac{t}{s+t}\beta_2\in K(\b_1,\b_2)$$
Thus $K(0,\alpha_i)K(\beta_1,\beta_2)=K(0,\alpha_i\beta_1,\alpha_i\beta_2)$ and
$K_1K_2=K(0,\alpha_2\beta_1,\alpha_2\beta_2)\cup K(0,\alpha_1\beta_1,\alpha_1\beta_2)$.

{\rm (2.b)} Let $x\in K(\a_1\b_1,\a_1\b_2,\a_2\b_1,\a_2\b_2)$. Then $x
=s\alpha_1\beta_1+t\alpha_1\beta_2+u\alpha_2\beta_1+v\alpha_2\beta_2
$ for some $  s,\ t,\ u,\ v\ge 0$ with $s+ t+ u+v=1$. Since $0 \notin K(\a_1,\a_2)$,
$\alpha_2=k\alpha_1$ for some $k>0$,
then $x=(p\alpha_1+(1-p)\alpha_2)(q\beta_1+(1-q)\beta_2),$ where $$p=s+t,  \quad
q=\frac{s+uk}{s+t+k(u+v)}\in [0,1]. $$
\vskip -.3in \qed

\medskip
The situation is more involved if neither
$K(0,\alpha_1,\alpha_2)$ nor $K(0,\beta_1,\beta_2)$ is a line segment. To describe the shape of $K_1K_2$ in such a case, we put the two segments
in a certain ``canonical'' position. More specifically, the next proposition shows that we can find $\alpha_0$ and $\beta_0\in\IC$ such that $\alpha_0^{-1}K_1$ and $\beta_0^{-1}K_2$ lie in the vertical line $\{z\in \IC: \Re(z)=1\}$.

\begin{proposition}\label{2.4}
Let $K_1 = K(\a_1, \a_2)$ and $K_2 = K(\b_1,\b_2)$ be two line segments in
$\IC$ such that neither $K(0,\a_1,\a_2)$ nor 
$K(0,\b_1,\b_2)$ is a line segment. Let
\begin{equation}\label{ab0}
\a_0=\dfrac{\a_1\overline{\a_2}-\a_2\overline{\a_1}}{2(\overline{\a_2}- \overline{\a_1})}\
\mbox{ and }\
\b_0=\dfrac{\b_1\overline{\b_2}-\b_2\overline{\b_1}}{2(\overline{\b_2}- \overline{\b_1})}
\end{equation}
Then $\a_0$ (respectively, $\b_0$) is the point on the line passing through $\a_1$ and $\a_2$
(respectively, $\b_1$ and $\b_2$) closest to 0. We have
\begin{equation}\label{ab1}
\dfrac{\a_1}{\a_0}=1+a_1i,\ \dfrac{\a_2}{\a_0}=1+a_2i, \
\dfrac{\b_1}{\b_0}=1+b_1i,\ \dfrac{\b_2}{\b_0}=1+b_2i
\end{equation}
for some $a_1,\ a_2,\ b_1$ and $b_2\in \IR$.
\end{proposition}

\it Proof. \rm  The line  passing through $\a_1$ and $\a_2$ is given by the parametric equation
 $r(t)= \a_1+t(\a_1-\a_2),\ t\in\IR$.  $\a_0$ in (\ref{ab0})
is obtained by minimizing $|r(t)|^2$.
Similarly, we have $\b_0$. By direct calculation we have   (\ref{ab1}) with
$$\begin{array}{rl}
a_1=&\dfrac{\a_1\overline{\a_2}+\a_2\overline{\a_1}-2|\a_1|^2}{i(\a_1\overline{\a_2}
-\a_2\overline{\a_1})},\
a_2=\dfrac{\a_1\overline{\a_2}+\a_2\overline{\a_1}-2|\a_2|^2}{i(\a_2\overline{\a_1}
-\a_1\overline{\a_2})},\\&\\
b_1=&\dfrac{\b_1\overline{\b_2}+\b_2\overline{\b_1}-2|\b_1|^2}{i(\b_1\overline{\b_2}
-\b_2\overline{\b_1})},\
b_2=\dfrac{\b_1\overline{\b_2}+\b_2\overline{\b_1}-2|\b_2|^2}{i(\b_2\overline{\b_1}
-\b_1\overline{\b_2})}.\end{array} $$ \vskip -0.3in
\qed

We can now describe  $K_1K_2$ for  two line
segments $K_1 =  K(\a_1,\a_2)$ and $K_2 =  K(\b_1,\b_2)$
in the  ``canonical'' position.
Because $K(\a_1,\a_2)K(\b_1,\b_2)$ is a simply connected
set, we focus on
the description of the boundary and the set of star centers of
 $K_1K_2$ in the following.

\begin{theorem} \label{thm2.4}
Let $K_1 = K(\a_1, \a_2)$ and $K_2 = K(\b_1,\b_2)$ with $\a_1 = 1+ia_1, \a_2 = 1+ia_2,
\b_1 = 1+ib_1, \b_2 = 1+ib_2$ such that $a_1 < a_2$
and $b_1 < b_2$. Assume $a_1 \le b_1$; otherwise,
interchange the roles of $K_1$ and $K_2$.
Then one of the following holds.

\begin{itemize}
\item[{\rm (a)}] $a_1 < a_2 \le b_1 < b_2$. Then
$K_1K_2$ is the convex  quadrilateral
$K(\a_1\b_1, \a_1\b_2,  \a_2\b_1, \a_2\b_2),$
which will degenerate to the triangle
$K(\a_1\b_1,a_1\b_2, \a_2\b_2)$
if $a_2 = b_1$; see Figure {\rm 1 (a) and (a.i)}.

\item
[{\rm (b)}] $a_1 \le b_1 <  a_2 \le b_2$.
Then $K_1K_2 \subseteq K(\a_1\b_1,a_1\b_2, \a_2\b_2)$, and  the boundary of
$K_1K_2$ consists of  the line segments
$K(\a_2^2, \a_2\b_2),$
$K(\a_2\b_2, \a_1\b_2),$
$K(\a_1\b_2, \a_1\b_1),$
$K(\a_1\b_1, \b_1^2),$
 and the curve 
$\mathbf{E}=\{(1+si)^2: s \in [b_1,a_2]\} \subseteq \bC.$
Here,
$K(\a_2^2, \a_2\b_2)$ lies on the tangent line of the curve
$\mathbf{E}$ at $\a_2^2$, and
$K(\b_1^2, \a_1\b_1)$ lies on the tangent line of the curve
$\mathbf{E}$ at $\b_1^2$.
The set of star centers equals
$K(\a_1,\b_1)K(\a_2,\b_2)$, which may be a quadrilateral,
a line or a point; see Figure
{\rm 1 (b), (b.i),
(b.ii), b(iii).}

\item[
{\rm (c)}] Suppose $a_1 < b_1 < b_2 < a_2$. Then
the boundary of
$K_1K_2$ consists of
the line segments
$K(\b_2^2, \a_2\b_2),$
$K(\a_2\b_2, \a_2\b_1),$
$K(\a_2\b_1, \b_1\b_2),$
$K(\b_1\b_2, \a_1\b_2),$
$K(\a_1\b_2,\a_1\b_1),$
$K(\b_1^2, \a_1\b_1)$
 and the curve segment
$\{(1+si)^2:s \in [b_1,b_2]\}\subseteq \bC.$
Here,
$K(\b_2^2, \a_2\b_2)$ lies on the tangent line of the curve
$\bC$ at $\b_2^2$, and
$K(\b_1^2, \a_1\b_1)$ lies on the tangent line of the curve
$\bC$ at $\b_1^2$. The unique star center is $\b_1\b_2$; see Figure {\rm 1 (c)}.

\end{itemize}
\end{theorem}

\begin{center}
\vspace{0.5cm}

\bigskip

%first figure b2=2 b1=1.5 a2=0.25 a1=-0.75
\begin{tabular}{cc}
\begin{tikzpicture}[line cap=round,line join=round,x=2.0cm,y=2.0cm]
\draw[opacity=0.5] (0.625,1.75)-- (0.5,2.25);
\draw[opacity=0.5] (2.5,1.25)-- (2.125,0.75);
\draw[opacity=0.5] (0.5,2.25)-- (2.5,1.25);
\draw[opacity=0.5] (0.625,1.75)-- (2.125,0.75);
\fill[opacity=0.2] (0.625,1.75)-- (0.5,2.25) -- (2.5,1.25) -- (2.125,0.75) -- cycle;
\begin{scriptsize}
\draw [fill=black] (0.5,2.25) circle (1.5pt);
\draw[color=black] (0.5,2.25) node[above] {$(1+b_2i)(1+a_2i)$};
\draw [fill=black] (2.5,1.25) circle (1.5pt);
\draw[color=black] (2.5,1.25) node[right] {$(1+b_2i)(1+a_1i)$};
\draw [fill=black] (0.625,1.75) circle (1.5pt);
\draw[color=black] (0.625,1.75) node[left] {$(1+b_1i)(1+a_2i)$};
\draw [fill=black] (2.125,0.75) circle (1.5pt);
\draw[color=black] (2.125,0.75) node[below] {$(1+b_1i)(1+a_1i)$};
\end{scriptsize}
\end{tikzpicture} \qquad & \quad
\begin{tikzpicture}[line cap=round,line join=round,x=1.5cm,y=1.5cm,scale=0.8]
\draw[opacity=0.5] (1.75,0.25)-- (0.5,2.25);
\draw[opacity=0.5] (-1,3)-- (0,2);
\draw[opacity=0.5] (0.5,2.25)-- (-1,3);
\draw[opacity=0.5] (1.75,0.25)-- (0,2);
\fill[opacity=0.2] (1.75,0.25)-- (0.5,2.25) -- (-1,3) -- (0,2) -- cycle;
\begin{scriptsize}
\draw [fill=black] (0.5,2.25) circle (1.5pt);
\draw[color=black] (0.5,2.25) node[right] {$(1+b_2i)(1+a_1i)$};
\draw [fill=black] (-1,3) circle (1.5pt);
\draw[color=black] (-1,3) node[above] {$(1+b_2i)(1+a_2i)$};
\draw [fill=black] (1.75,0.25) circle (1.5pt);
\draw[color=black] (1.75,0.25) node[left] {$(1+a_2i)(1+a_1i)$};
\draw [fill=black] (0,2) circle (1.5pt);
\draw[color=black] (0,2) node[anchor= north east] {$(1+a_2i)^2$};
\end{scriptsize}
\end{tikzpicture}
\\

\textbf{(a)} $a_1<a_2<b_1<b_2$, & \textbf{(a.i)} $a_1<a_2=b_1 < b_2$.\\

\\

\begin{tikzpicture}[line cap=round,line join=round,x=1.3cm,y=1.3cm]
\begin{tiny}
\fill (1,2) circle (2pt);
\draw (1,2) node[right] {$(1+a_1i)(1+b_2i)$};
\fill (1,1.6) circle (2pt);
\draw (1,1.6) node[right] {$(1+a_1i)(1+a_2i)$};
\fill (1,0.5) circle (2pt);
\draw (1,0.5) node[right] {$(1+a_1i)(1+b_1i)$};
\draw (0,2.5) node[anchor= south west] {$(1+b_1i)(1+b_2i)$};
\fill (0,2.5) circle (2pt);
\draw (-2.2,3.6) node[above] {$(1+a_2i)(1+b_2i)$};
\fill (-2.2,3.6) circle (2pt);
\draw (-1.56,3.2) node[left] {$(1+a_2i)^2$};
\fill (-1.56,3.2) circle (2pt);
\fill (0.75,1) circle (2pt);
\draw (0.75,1) node[left]  {$(1+b_1i)^2$};
\fill (0.2,2.1) circle (2pt);
\draw (-0.2,2) node[left] {$(1+a_2i)(1+b_1i)i$};
\end{tiny}
\draw[scale=1,domain=1:3.2,smooth,variable=\y,opacity=0.5]  plot ({1-(\y*\y)/4},{\y});
\draw[opacity=0.5] (-2.2,3.6) -- (1,1.6);
\draw[opacity=0.5] (1,2) -- (1,0.5);
\draw[opacity=0.5] (1,2) -- (-2.2,3.6);
\draw[opacity=0.5] (1,0.5) -- (0,2.5);
\fill[scale=1,domain=1:3.2,smooth,variable=\y,line width=1.2,opacity=0.2]
plot ({1-(\y*\y)/4},{\y}) -- (-2.2,3.6) --  (1,2) -- (1,0.5) --  cycle;
\draw[->,help lines,shorten >=4pt] (-0.2,2) .. controls (-0.1,1.9) and (0,2.1) .. (0.2,2.1);

\end{tikzpicture} &

\begin{tikzpicture}[line cap=round,line join=round,x=1.0cm,y=1.0cm,scale=0.8]
\draw[smooth,samples=100,domain=-2.0:2.0,opacity=0.5] plot ({1.0-0.25*\x*\x},{\x});
\draw[opacity=0.5] (-1.,3.)-- (3.,1.);
\draw[opacity=0.5] (0.,-2.)-- (2.,0.);
\draw[opacity=0.5] (3.,1.)-- (0.,-2.);
\draw[opacity=0.5] (-1.,3.)-- (2.,0.);
\fill[opacity=0.2,smooth,samples=100,domain=-2.0:2.0] plot ({1.0-0.25*\x*\x},{\x}) -- (-1.,3.)--
(3.,1.);
\begin{scriptsize}
\draw [fill=black] (0.,-2.) circle (1.5pt);
\draw[color=black] (0.,-2.) node[left] {$(1+a_1i)^2$};
\draw [fill=black] (0.,2.) circle (1.5pt);
\draw[color=black] (0.,2.) node[anchor = north east] {$(1+a_2i)^2$};
\draw [fill=black] (2.,0.) circle (1.5pt);
\draw[color=black] (2.,0.) node[anchor= north west] {$(1+a_1i)(1+a_2i)$};
\draw [fill=black] (-1.,3.) circle (1.5pt);
\draw[color=black] (-1.,3.) node[right] {$(1+a_2i)(1+b_2i)$};
\draw [fill=black] (3.,1.) circle (1.5pt);
\draw[color=black] (3.,1.) node[right] {$(1+a_1i)(1+b_2i)$};
\end{scriptsize}
\end{tikzpicture} \\

\textbf{(b)} $\ a_1 <  b_1 < a_2 < b_2$,  & \textbf{(b.i)} $a_1=b_1<a_2<b_2$, \\

\\

 \begin{tikzpicture}[line cap=round,line join=round,x=0.5cm,y=0.75cm]
\draw[smooth,samples=100,domain=0.0:4.0,opacity=0.5] plot({1.0-0.25*\x*\x},{\x});
\draw[opacity=0.5] (-3.,4.)-- (1.,2.);
\draw[opacity=0.5] (1.,-1.)-- (3.,1.);
\draw[opacity=0.5] (1.,2.)-- (1.,-1.);
\draw[opacity=0.5] (-3.,4.)-- (3.,1.);
\fill[smooth,samples=100,domain=0.0:4.0,opacity=0.2] plot({1.0-0.25*\x*\x},{\x}) -- (3,1) --
(1,-1) -- cycle;
\begin{scriptsize}
\draw [fill=black] (1.,0.) circle (1.5pt);
\draw[color=black] (1.,0.) node[left] {$(1+b_1i)^2$};
\draw [fill=black] (-3.,4.) circle (1.5pt);
\draw[color=black] (-3.,4.) node[right] {$(1+a_2i)^2$};
\draw [fill=black] (1.,2.) circle (1.5pt);
\draw[color=black] (1.,2.) node[anchor= south west] {$(1+a_2i)(1+b_1i)$};
\draw [fill=black] (1.,-1.) circle (1.5pt);
\draw[color=black] (1.,-1.) node[left] {$(1+a_1i)(1+b_1i)$};
\draw [fill=black] (3.,1.) circle (1.5pt);
\draw[color=black] (3.,1.) node[right] {$(1+a_1i)(1+a_2i)$};
\end{scriptsize}
\end{tikzpicture}
 &
%here a2=1 and a1=-0.75
\begin{tikzpicture}[line cap=round,line join=round,x=1.0cm,y=1.0cm]
\draw[opacity=0.5] (0.,2.)-- (1.75,0.25);
\draw[opacity=0.5] (0.4375,-1.5)-- (1.75,0.25);
\draw[smooth,samples=100,domain=-1.5:2.0,opacity=0.5] plot ({1-0.25*\x*\x},{\x});
\fill[smooth,samples=100,domain=-1.5:2.0,opacity=0.2] plot ({1-0.25*\x*\x},{\x})
-- (1.75,0.25) --  cycle;
\begin{scriptsize}
\draw [fill=black] (0.,2.) circle (1.5pt);
\draw[color=black] (0,2) node[left] {$(1+a_2i)^2$};
\draw [fill=black] (0.4375,-1.5) circle (1.5pt);
\draw[color=black] (0.4375,-1.5) node[left] {$(1+a_1i)^2$};
\draw [fill=black] (1.75,0.25) circle (1.5pt);
\draw[color=black] (1.75,0.25) node[right] {$(1+a_1i)(1+a_2i)$};
\end{scriptsize}
\end{tikzpicture}

\\

\textbf{(b.ii)} $a_1<b_1<a_2=b_2$, &
\textbf{(b.iii)} $a_1 = b_1<b_2 = a_2$. \\
\end{tabular}

\bigskip
\begin{tabular}{c}
\begin{tikzpicture}[line cap=round,line join=round,x=1.0cm,y=1.0cm]
\draw[smooth,samples=100,domain=0.0:2.0,opacity=0.5] plot({1-0.25*\x*\x},{\x});
\draw[opacity=0.5] (-1.,3.)-- (1.,2.);
\draw[opacity=0.5] (1.,-1.)-- (2.,0.);
\draw[opacity=0.5] (1.,2.)-- (1.,-1.);
\draw[opacity=0.5] (-1.,3.)-- (2.,0.);
\fill[opacity=0.2,smooth,samples=100,domain=0.0:2.0]  plot({1-0.25*\x*\x},{\x}) -- (-1.,3.) --
(1,2) -- (1,1) -- (2,0) -- (1,-1) -- cycle;
\begin{scriptsize}
\draw [fill=black] (1,0.) circle (1.5pt);
\draw[color=black] (1,0) node[left] {$(1+b_1i)^2$};
\draw [fill=black] (0,2) circle (1.5pt);
\draw[color=black] (0,2) node[anchor=north east] {$(1+b_2i)^2$};
\draw [fill=black] (1,1) circle (1.5pt);
\draw[color=black] (1,1) node[right] {$(1+b_1i)(1+b_2i)$};
\draw [fill=black] (-1,3) circle (1.5pt);
\draw[color=black] (-1,3) node[above] {$(1+a_2i)(1+b_2i)$};
\draw [fill=black] (1,-1) circle (1.5pt);
\draw[color=black] (1,-1) node[left] {$(1+a_1i)(1+b_1i)$};
\draw [fill=black] (2,0) circle (1.5pt);
\draw[color=black] (2,0) node[right] {$(1+a_1i)(1+b_2i)$};
\draw [fill=black] (1,2) circle (1.5pt);
\draw[color=black] (1,2) node[right] {$(1+a_2i)(1+b_1i)$};
\end{scriptsize}
\end{tikzpicture}  \\

\textbf{(c)} $a_1 < b_1<b_2 < a_2$. \\

\end{tabular}

\medskip
\textbf{Figure 1.} The set $K(1+a_1i, 1+a_2i)K(1+b_1i, 1+b_2i)$
described in Theorem \ref{thm2.4}.

\end{center}

\bigskip
To prove Theorem \ref{thm2.4}, we need the following lemma that treat some special
cases of the theorem. It turns out that these special cases are the building blocks
for the general case.

\medskip
\begin{lemma}\label{lem2.3} Let $a_1<a_2\leq b_1<b_2$. Then
\begin{itemize}
\item[{\rm (a)}] $K(1+a_1i,1+a_2i)K(1+b_1i,1+b_2i)$ is
the quadrilateral (or triangle if $a_2=b_1$),
$$ \mathbf{K}=K\Big((1+a_1i)(1+b_1i),(1+a_1i)(1+b_2i),(1+a_2i)(1+b_1i),(1+a_2i)(1+b_2i)\Big).$$
\item[{\rm (b)}] $K(1+a_1i,1+a_2i)K(1+a_1i,1+a_2i)$ is the simply connected region
bounded by
the line segments $$\mathbf{L}_1=K\Big((1+a_1i)^2,(1+a_1i)(1+a_2i)\Big),\quad \mathbf{L}_2
=K\Big((1+a_2i)^2,(1+a_1i)(1+a_2i)\Big),$$ 
and the curve
$\mathbf{E}=\{(1+si)^2\ : \ s\in [a_1,a_2]\}$. The set $\mathbf{L}_1$ is a segment of
the tangent line of $\mathbf{E}$ at $(1+a_1i)^2$, and $\mathbf{L}_2$ is a segment of the
tangent line of $\mathbf{E}$ at $(1+a_2i)^2$.
\end{itemize}
\end{lemma}

 \medskip\noindent
{\it Proof.}
(a) Suppose $\alpha_j=1+a_ji$ and $\beta_j=1+b_ji$ for $j=1,2$
are such that $a_1 < a_2 \le b_1 < b_2$. Let $K_1=K(\alpha_1,\alpha_2)$ and $K_2=K(\beta_1,\beta_2)$.
It suffices to show that the union of the line segments
$$\ell_1=\beta_2  K_1,\quad \ell_2=\beta_1 K_1,\quad \ell_3=\alpha_2  K_2,
\quad  \ell_4=\alpha_1 K_2$$
forms the boundary of the  quadrilateral (or triangle) $\mathbf{K}$, that is, the union is a simple closed curve. By simply connectedness and the fact that $ K_1K_2$ is a subset of $\mathbf{K}$, we get the desired conclusion. For the convenience of discussion, we will identify $x+iy\in\IC$ with $(x,y)\in \IR^2$ and $(x,y,0)\in\IR^3$.
Note that since $\arg(\alpha_1\beta_1)< \arg(\alpha_2\beta_1), \arg(\alpha_1\beta_2)< \arg(\alpha_2\beta_2)$, it suffices to show that $\alpha_1\beta_2$ and $\alpha_2\beta_1$ are on opposite sides of the line  $\ell$ passing through $\alpha_1\beta_1$ and $\alpha_2\beta_2$. This is true if and only if the cross product $ (\a_2\b_1-\a_2\b_2) \times   (\a_1\b_1-\a_2\b_2)$ and $(\a_1\b_2-\a_2\b_2)\times  (\a_1\b_1-\a_2\b_2)$ are pointing in opposite directions, that is
%the cross products $\ell_1\times \ell$ and $\ell_2\times \ell$ are pointing in opposite directions, that is
$$\det\begin{bmatrix}
\Re(\alpha_2\beta_1-\alpha_2\beta_2) & \Re(\alpha_1\beta_1-\alpha_2\beta_2)\\
\Im(\alpha_2\beta_1-\alpha_2\beta_2) & \Im(\alpha_1\beta_1-\alpha_2\beta_2)
\end{bmatrix} \cdot \det\begin{bmatrix}
\Re(\alpha_1\beta_2-\alpha_2\beta_2) & \Re(\alpha_1\beta_1-\alpha_2\beta_2)\\
\Im(\alpha_1\beta_2-\alpha_2\beta_2) & \Im(\alpha_1\beta_1-\alpha_2\beta_2)
\end{bmatrix} \leq 0$$
The expression on the left hand side is
$$ [(b_1-b_2)(a_2-a_1)(a_2-b_1)]\cdot[(b_1-b_2)(a_2-a_1)(b_2-a_1)]=(b_1-b_2)^2(a_2-a_1)^2(a_2-b_1)(b_2-a_1)$$
Since $a_2\leq b_1$ and $b_2>a_1$, then we are done.

\medskip\noindent
\medskip\noindent
To prove (b), first note that $\mathbf{L}_1,\mathbf{L}_2$ and $\mathbf{E}$ are clearly in $K_1K_1$. Direct calculation shows that $\mathbf{L}_1$ with equation $x= 1-a_1(y-a_1)$  and $\mathbf{L}_2$ with equation $x=1-a_2(y-a_2)$ are tangent to the parabola $\mathbf{E}$ with equation $x=1-\frac{y^2}{4} $ at the points $(1-a_1^2,\ 2a_1)$ and $(1-a_2^2,\ 2a_2)$ respectively.
\\
Since $K_1K_1$ is simply connected, the region
\begin{equation}\label{kset}\mathbf{S}=\left\{x+iy:1-\frac{y^2}{4}\le  x\le 1-a_1(y-a_1),\ 1-a_2(y-a_2)
  \right\},
  \end{equation}
which is the region enclosed by $\mathbf{L}_1,\mathbf{L} _2$ and $\mathbf{E}$ is a subset of $K_1K_1$. Now, suppose  $x+iy\in K_1K_1$. Then there exist $r$ and $s$ with $a_1\leq r,s\leq a_2$ such that $$x+iy=(1+ir)(1+is)=1-rs+i(r+s).$$   Note that
$$x=1-rs\ge 1-\dfrac{1}{4}(r+s)^2=1-\dfrac{y^2}{4}$$
always holds. Also, if $a\le t\le b$, then   $(a+b-t)t\ge ab$.
Since
$$ a_1\leq r \leq s+r-a_1\mbox{ and }s+r-a_2\leq r \leq  a_2 \,,$$
we have $rs\ge a_1(s+r-a_1),\ a_2(s+r-a_2)$. Hence,
$$\begin{array}{rl}x=&1-rs\le 1-a_1(r+s-a_1)=1-a_1(y-a_1),\quad \mbox{ and }  \\&\\ x=&1-rs\le 1-a_2(r+s-a_2)=1-a_2(y-a_2).\end{array}$$
This shows that $K_1K_1$ lies inside $\mathbf{S}$. Thus $K_1K_1=\mathbf{S}$. 
 \qed

\it Proof of Theorem \ref{thm2.4}. \rm  Suppose
$K_1 = K(1+ia_1, 1+ia_2)$ and $K_2 = K(1+ib_1, a+ib_2)$ such that
$a_1 \le a_2, b_1 \le b_2$. We show that if $K_1K_2$ can be written as
the union of subsets of the form in Lemma \ref{lem2.3}. In fact,
if $[a_1,a_2]\cap [b_1,b_2]=[c_1,c_2]$, then
$$ K_1K_2=(\a_0\b_0)\left[(A  C)\cup (A  B) \cup
(C  C) \cup (C  B)\right],$$
where $C=K(1+c_1i,1+c_2i)$, $B=K(1+b_1i,1+b_2i)\setminus C$ and $A=K(1+a_1i,1+a_2i)
\setminus C$. By Lemma \ref{lem2.3}, we get the conclusion.
\qed

By Theorem \ref{thm2.4}, we have the following corollary giving information about the
star center of the product of two line segments without putting them in
the  ``canonical'' position.

\begin{corollary}\label{2.6} Let $K_1=K(\a_1,\a_2)$ an $K_2=K(\b_1,\b_2)$,
where $\a_1,\a_2,\b_1,\b_2\in \mathbb{C}$ such that
$\arg(\alpha_1)<\arg(\alpha_2)< \arg(\alpha_1)+ \pi$ and $\arg(\beta_1)<\arg(\beta_2)
< \arg(\beta_1)+ \pi$.
Then $K_1K_2$ is star-shaped and one of the following holds.
\begin{enumerate}
\item [{\rm (a)}]
There exists $\xi \in \mathbb{C}$ such that $\xi K_1\subseteq K_2$. Equivalently,
the segments
$K(\alpha_1\beta_1,\alpha_1\beta_2)$ and $K(\alpha_2\beta_1,\alpha_2\beta_2)$
intersect at  $\xi\alpha_1\alpha_2$.
In this case, $\xi\alpha_1\alpha_2$ is the unique star-center  of $K_1K_2$.
\item [{\rm (b)}]There exists $\xi \in \mathbb{C}$ such that $\xi K_2\subseteq K_1$. Equivalently, the
segments $K(\alpha_1\beta_1,\alpha_2\beta_1)$ and $K(\alpha_1\beta_2,\alpha_2\beta_2)$
intersect  at $\xi\beta_1\beta_2$.
In this case, $\xi\beta_1\beta_2$ is the unique star-center of $K_1K_2$.
\item[{\rm (c)}] Condition {\rm (a)} and {\rm (b)}
do not hold, and  every point in $K(\beta_1\alpha_2,\beta_2\alpha_1)$ is a star center of $K_1K_2$
\end{enumerate}
\end{corollary}

\section{The product set of two convex polygons}

In this section, we study the product set of two convex polygons (including interior).
It is known that  for every convex polygon $K_1$ with
vertexes $\mu_1, \dots, \mu_n$, then
$K_1 = W(T)$ for $T = \diag(\mu_1, \dots, \mu_n) \in M_n$.
In Section 3.1, we will show that
the product set of two convex polygons may not be star-shaped.
In particular, we have a product set of two triangles that are not star-shaped.
This gives a negative answer to
the conjecture in \cite{Karol}.

\subsection{Products of polygons that are not star-shaped}

In this subsection, we show that there are examples  $K_1$ and $K_2$
such that $K_1K_2$ is not star-shaped. The first example
has the form $K_1 = K_2 = K(\a_1, \bar \a_1, \a_2)$, where $\a_2\notin \mathbf{R}$.
One can regard $K_1 = W(T)$ with
$T = \diag(\a_1, \bar \a_1, \a_2) \in M_3$ so that
 the set
$W^{\otimes}(T\otimes T) = W(T)W(T)$ is not star-shaped.
We can construct another example of the form $K_1 = K_2 = K(\a_1, \bar \a_1, \a_2,\bar \a_2)$,
which is symmetric about the real axis, such that $K_1K_2$ is not star-shaped. One can regard
$K_1 = W(A)$ for a real normal
matrix $A \in M_4$ with eigenvalues $\a_1, \bar \a_1, \a_2,\bar \a_2$
so that  $W^{\otimes}(A\otimes A)$ is not star-shaped.

\begin{example} \label{3.1} Let $K_1=K(e^{i\frac{\pi}{3}},e^{-i\frac{\pi}{3}},
0.95e^{i\frac{\pi}{4}})$.
Then $K_1K_1$ is not star-shaped.

\it Proof. \rm
Let $\a_1=e^{i\frac{\pi}{3}}$ and $\a_2=0.95e^{i\frac{\pi}{4}}$,
$K_1= K(\a_1,\overline{\a_1}, \a_2)$. Then
$1=\a_1 \overline{\a_1},\ 0.95^2i=\a_2^2\in K_1K_1$.
We are going to show that {\bf  a)} if $s$ is a star center of $K_1K_1$, then $s=1$ and {\bf b)} $(1-t)+t0.95^2i\not \in K_1K_1$ for all $t\in (0,1)$.

Let $S$ be a closed and bounded subset of $\IC$, with $0\not\in S$. Suppose $t\in \IR$ and
$S\cap\{re^{it}:r>0\}\not = \emptyset$. Let $\rho_0^S(t)=\min\{r>0:re^{it}\in S\}$ and
$\rho_1^S(t)=\max\{r>0:re^{it}\in S\}$.

Let $L_1= K(\a_1,\overline{\a_1})$, $S_1=K_1K_1$ and $S_2=L_1L_1$. Since $\rho_0^{K_1}(\t)=\rho_0^{L_1}(\t)$  for $-\dfrac{\pi}{3} \le \t\le \dfrac{\pi}{3} $,   it follows that $\rho_0^{S_1}(\t)= \rho_0^{S_2}(\t) $ for $-\dfrac{2\pi}{3} \le \t\le \dfrac{2\pi}{3} $.

Note that $x+iy\in S_2\Lra 4(x+iy)\in (2L_1)(2L_1)$. Then, applying Lemma \ref{lem2.3} (b) to $2L_1=K(1-i\sqrt{3},1+i\sqrt{3})$, we have
$$S_2=\{x+iy:1-4y^2\le 4x\le 1-\sqrt{3}(4y-\sqrt{3}),\ 1+\sqrt{3}(4y+\sqrt{3})\}$$

\begin{center}
\begin{tikzpicture}[line cap=round,line join=round,x=1.0cm,y=1.0cm,scale=2.5]
\draw [samples=50,variable=\y,domain=-0.8660254037844386:0.8660254037844386] plot ({1/4-\y*\y},\y);
\draw (-0.1,0) node {\scriptsize $4x=1-4y^2$};
\draw (-0.5,0.8660254037844386)-- (1.5,-0.29);
\draw (-0.5,0.8660254037844386)-- (1,0) node[midway,above, sloped] {\scriptsize $4x=1-\sqrt{3}(4y-\sqrt{3})$};
\draw (-0.5,-0.8660254037844386)-- (1,0) node[midway,below, sloped] {\scriptsize $4x=1+\sqrt{3}(4y+\sqrt{3})$};
\draw (1.5,0.29)-- (-0.5,-0.8660254037844386);
\fill[fill opacity=0.2,samples=50,variable=\y,domain=-0.8660254037844386:0.8660254037844386] plot ({1/4-\y*\y},\y) -- (-0.5,0.8660254037844386)-- (1.,0.) -- cycle;
\fill[pattern=north west lines] (1,0) -- (1.5,0.29) -- (1.5,-0.29) --cycle;
\begin{footnotesize}
\fill  (-0.5,0.8660254037844386) circle (1pt);
\draw (-0.5,0.8660254037844386) node[above] {$e^{i\frac{2\pi}{3}}$};
\fill  (1,0.) circle (1pt);
\draw (1,-0.1) node[below] {$1$};
\fill  (-0.5,-0.8660254037844386) circle (1pt);
\draw (-0.5,-0.8660254037844386) node[below] {$e^{-i\frac{2\pi}{3}}$};
\fill[fill=white]  (1.3,0) circle (3pt);
\draw (1.3,0) node {$s$};
\draw (0.5,0) node {$S_2$};
\end{footnotesize}
\end{tikzpicture}\\
\textbf{Figure 2.} Plot of $S_2=L_1L_1$
\end{center}

{\bf  a)} Note that $\Big\{\rho_0^{S_1}(\theta)\ :\ \theta\in [-2\pi/3, 2\pi/3]\Big\} =\Big\{\rho_0^{S_2}(\theta)\ :\ \theta\in [-2\pi/3, 2\pi/3]\Big\}=
\{z^2\ : \ z\in L_1\}$. This means that the curve $\{z^2\ :\ z\in L_1\}$ is a boundary curve of $S_2$. By Proposition 1.2, if $s$ were a star-center of $S_2$, then the segment $K(s,z^2)$ must be in $S_2$ for any $z\in L_1$.

If $s=x+iy$ is a star center of $S_1$, then we must have
$$4x\ge 1-\sqrt{3}(4y-\sqrt{3}),\ 1+\sqrt{3}(4y+\sqrt{3})\Ra x\ge 1$$
Since $|z|\le 1$ for all $z\in S_1$, we have $s=1$.

{\bf b)} Let  $L_2 = K(\a_1, \a_2)$, $L_3 = K(\bar \a_1, \a_2)$.
Then the boundary of the simply connected set
$S_1=K_1K_1$ is a subset of $\cup_{1 \le i\le j \le 3} L_i L_j$.

Suppose $0< \t< \frac{\pi}{2}$ and $\rho_1^{S_1}(t)=r$.
Then $re^{i\t}\in L_2L_3\cup L_3L_3$. Direct calculation using Lemma \ref{lem2.3}
and Proposition \ref{2.4} shows that $\rho_1^{L_2L_3}(\t),\ \rho_1^{L_3L_3}(\t)
<\rho_1^{K(1,\a_2^2)}(\t)$; see the following figures.

\begin{center}
\begin{tabular}{ccc}
\begin{tikzpicture}[line cap=round,line join=round,x=1.0cm,y=1.0cm,scale=2.5]
\draw[->,color=black] (-0.5,0.) -- (1.25,0.);
\draw[->,color=black] (0.,-1) -- (0.,1.25);
\draw  (-0.2458780928473946,0.9176295349746149) -- (0.,0.9025);
\draw (0.,0.9025)-- (0.9176295349746149,-0.2458780928473947);
\draw  (-0.2458780928473946,0.9176295349746149) -- (1.,0.);
\draw (1.,0.)-- (0.9176295349746149,-0.2458780928473947);
\fill[fill opacity=0.2]  (-0.2458780928473946,0.9176295349746149) -- (0.,0.9025)--(0.9176295349746149,-0.2458780928473947) -- (1,0) -- cycle;
\draw[samples=50,variable=\y,domain=0.1230:0.67] plot ({0.46-0.3234*\y-0.46*(\y*\y)/4},{0.3234+0.46*\y-0.3234*(\y*\y)/4});
\fill[fill opacity=0.2,samples=50,variable=\y,domain=0.1230:0.67] 
(0.3223115666172949,0.49913946280361987) -- 
plot ({0.46-0.3234*\y-0.46*(\y*\y)/4},{0.3234+0.46*\y-0.3234*(\y*\y)/4}) -- cycle;
\draw[dashed] (0,0.9025)--(1,0);
\begin{scriptsize}
\fill (1.,0.) circle (1pt);
\draw (1.1,0) node[above] {$1$};
\fill (0.,0.9025) circle (1pt);
\draw (0.,1) node[right] {$\alpha_2^2$};
\fill  (-0.2458780928473946,0.9176295349746149)  circle (1pt);
\draw  (-0.3,0.9176295349746149)  node[above] {$\alpha_1\alpha_2$};
\fill (0.9176295349746149,-0.2458780928473947) circle (1pt);
\draw (0.9176295349746149,-0.3)node[below] {$\alpha_2\overline{\alpha_1}$};
\end{scriptsize}
\end{tikzpicture}
 & \begin{tikzpicture}[line cap=round,line join=round,x=1.0cm,y=1.0cm,scale=2.5]
\draw[->,color=black] (-0.8,0.) -- (1.25,0.);
\draw[->,color=black] (0.,-1) -- (0.,1.25);
\draw (0.,0.9025)-- (0.9176295349746149,-0.2458780928473947);
\draw (0.9176295349746149,-0.2458780928473947)-- (-0.5,-0.8660254037844386);
\draw[samples=50,variable=\y,domain=-2.7154:2.5029] plot ({0.343+0.07761*\y-0.343*(\y*\y)/4},{-0.07761+0.3430*\y+0.07761*(\y*\y)/4});
\fill[fill opacity=0.2,samples=50,variable=\y,domain=-2.7154:2.5029] plot ({0.343+0.07761*\y-0.343*(\y*\y)/4},{-0.07761+0.3430*\y+0.07761*(\y*\y)/4}) --(0.9176295349746149,-0.2458780928473947)-- cycle ;
\draw[dashed] (0,0.9025)--(1,0);
\begin{scriptsize}
\fill (0.,0.9025) circle (1pt);
\draw (0.,0.9025) node[anchor=south west] {$\alpha_2^2$};
\fill (0.9176295349746149,-0.2458780928473947) circle (1pt);
\draw (0.9176295349746149,-0.2458780928473947) node[right] {$\alpha_2\overline{\alpha_1}$};
\fill (-0.5,-0.8660254037844386) circle (1pt);
\draw (-0.5,-0.8660254037844386) node[below] {$\overline{\alpha_1}^2$};
\end{scriptsize}
\end{tikzpicture} & \begin{tikzpicture}[line cap=round,line join=round,x=1.0cm,y=1.0cm,scale=2.5]
\draw[->,color=black] (-0.8,0.) -- (1.25,0.);
\draw[->,color=black] (0.,-1) -- (0.,1.25);

\draw (0.,0.9025)-- (0.9176295349746149,-0.2458780928473947);
\draw (0.9176295349746149,-0.2458780928473947)-- (-0.5,-0.8660254037844386);
\draw[samples=50,variable=\y,domain=-2.7154:2.5029] plot ({0.343+0.07761*\y-0.343*(\y*\y)/4},{-0.07761+0.3430*\y+0.07761*(\y*\y)/4});
\fill[fill opacity=0.2,samples=50,variable=\y,domain=-2.7154:2.5029] plot ({0.343+0.07761*\y-0.343*(\y*\y)/4},{-0.07761+0.3430*\y+0.07761*(\y*\y)/4}) --(0.9176295349746149,-0.2458780928473947)-- cycle ;

\draw  (-0.2458780928473946,0.9176295349746149) -- (0.,0.9025);
\draw  (-0.2458780928473946,0.9176295349746149) -- (1.,0.);
\draw (1.,0.)-- (0.9176295349746149,-0.2458780928473947);
\fill[fill opacity=0.2]  (-0.2458780928473946,0.9176295349746149) -- (0.,0.9025)--(0.9176295349746149,-0.2458780928473947) -- (1,0) -- cycle;
\draw[samples=50,variable=\y,domain=0.1230:0.67] plot ({0.46-0.3234*\y-0.46*(\y*\y)/4},{0.3234+0.46*\y-0.3234*(\y*\y)/4});
\fill[fill opacity=0.2,samples=50,variable=\y,domain=0.1230:0.67] 
(0.3223115666172949,0.49913946280361987) -- 
plot ({0.46-0.3234*\y-0.46*(\y*\y)/4},{0.3234+0.46*\y-0.3234*(\y*\y)/4}) -- cycle;

\draw [samples=50,variable=\y,domain=-0.8660254037844386:0.8660254037844386] plot ({1/4-\y*\y},\y);
\draw (-0.5,0.8660254037844386)-- (1,0);
\draw (1,0)-- (-0.5,-0.8660254037844386);
\fill[fill opacity=0.2,samples=50,variable=\y,domain=-0.8660254037844386:0.8660254037844386] plot ({1/4-\y*\y},\y) -- (-0.5,0.8660254037844386)-- (1.,0.) -- cycle;

\draw (-0.5,0.8660254037844386) --(-0.2458780928473946,0.9176295349746149)  ;
\draw (-0.2458780928473946,0.9176295349746149)   -- (0.9176295349746149,-0.2458780928473947);

\fill[fill opacity=0.2] (-0.5,0.8660254037844386) --(-0.2458780928473946,0.9176295349746149)   -- (0.9176295349746149,-0.2458780928473947) -- (1,0) -- cycle;

\fill[fill opacity=0.2] (-0.2458780928473946,0.9176295349746149) -- (1,0) -- (0.9176295349746149,-0.2458780928473947) --cycle;

\draw[dashed] (0,0.9025)--(1,0);
\begin{scriptsize}
\fill (0.,0.9025) circle (1pt);
\draw (0.,0.9025) node[anchor=south west] {$\alpha_2^2$};
\fill (0.9176295349746149,-0.2458780928473947) circle (1pt);
\draw (0.9176295349746149,-0.2458780928473947) node[right] {$\alpha_2\overline{\alpha_1}$};
\fill (-0.5,-0.8660254037844386) circle (1pt);
\draw (-0.5,-0.8660254037844386) node[below] {$\overline{\alpha_1}^2$};
\fill (1.,0.) circle (1pt);
\draw (1.1,0) node[above] {$1$};
\fill  (-0.2458780928473946,0.9176295349746149)  circle (1pt);
\draw  (-0.3,0.9176295349746149)  node[above] {$\alpha_1\alpha_2$};
\fill  (-0.5,0.8660254037844386) circle (1pt);
\draw (-0.5,0.8660254037844386) node[anchor= south east] {$\alpha_1^2$};
\end{scriptsize}
\end{tikzpicture}\\
Plot of $L_2L_3$ &  Plot of $L_3L_3$ & Plot of $K_1K_1$
\end{tabular}
\smallskip\\
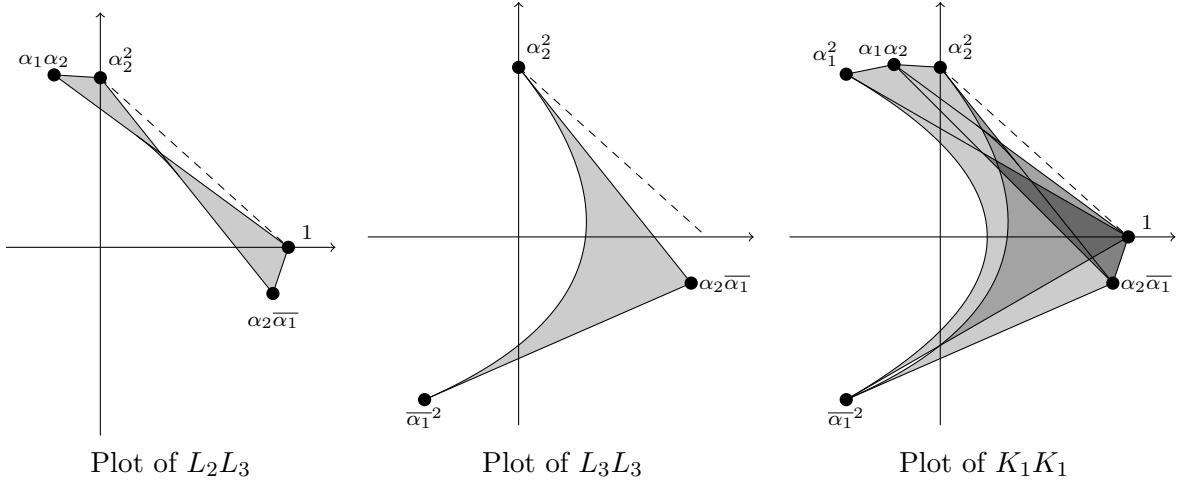
\textbf{Figure 3.}
\end{center}
\medskip
We conclude that $K_1K_1$ is not star-shaped. \qed
\end{example}

Next, we modify Example \ref{3.1} to Example \ref{3.2} so that
$\bar K_1 = K_1 (\a_1,\a_2,\bar \a_1,\bar \a_2)$ with
$\a_1=e^{i\frac{\pi}{3}}$ and $\a_2=0.95e^{i\frac{\pi}{4}}$.
In this case, one can regard $K_1=W(A)$ for some real symmetric $A\in M_4$. The product
set $K_1K_2$ will be larger than the product set considered in Example \ref{3.1}.
Never-the-less, we
can analyze the product of the sets $L_iL_j$ for $i,j=1,2,3,4$, where
$L_1=K(\a_1,\bar{\a}_1)$,  $L_2=K(\a_1,\a_2)$, $L_3=K(\a_2,\bar{\a}_2)$,
$L_4=K(\bar{\a}_2,\bar{\a}_1)$ so that $\cup_{1 \le i \le j \le 4} L_iL_j$
contains the boundary of the simply connected set $K_1K_1$.
Again one can show that the part of the boundary
$\{z^2: z \in K(\a_1,\bar\a_1)\}$ of $L_1L_1$ is also part of the boundary of $K_1K_1$
so that  $1 = \a_1\bar{\a}_1 \in K_1K_1$
is the only possible candidate to serve as a star-center for $K_1K_1$.
However,  none of the set $L_iL_j$ contains the set $\{t + (1-t)0.95^2i: 0 < t < 1/3\}.$
Thus, the line segment joining $1$ and $\a_2^2 = 0.95^2i$
is not in $K_1K_1$. Hence, $1$ is not the star center of $K_1K_1$,
and $K_1K_1$ is not star-shaped.

\begin{example} \label{3.2}
\rm
Let $K_1=K(e^{i\frac{\pi}{3}},e^{-i\frac{\pi}{3}},
0.95e^{i\frac{\pi}{4}},0.95e^{-i\frac{\pi}{4}})$.
Then $K_1$ is symmetric about the $x$-axis but $P=K_1K_1$ is not star-shaped.
\end{example}

\iffalse
\it Proof. \rm
Let $S$ be a closed and bounded subset of $\IC$, with $0\not\in S$. Suppose $t\in \IR$ and
$S\cap\{re^{it}:r>0\}\not = \emptyset$. Let $\rho_0^S(t)=\min\{r>0:re^{it}\in S\}$ and
$\rho_1^S(t)=\max\{r>0:re^{it}\in S\}$.

Suppose $K_2=K(e^{i\frac{\pi}{3}},e^{-i\frac{\pi}{3}})$ and $P_0=K_2K_2$.
Then by Lemma \ref{lem2.3} and Proposition \ref{2.4}, $1$  is the unique star center of $P_0$.

Let $P=K_1K_1$. Then for $-\dfrac{2\pi}{3}\le t\le \dfrac{2\pi}{3}$, we have
$$\begin{array}{rl}\rho_0^P(t)=&
\min\{\rho_0^{K_1}(t-s)\rho_0^{K_1}(s):-\dfrac{\pi}{3}\le t-s,\ s\le \dfrac{\pi}{3}\}
\\&\\=&\min\{\rho_0^{K_2}(t-s)\rho_0^{K_2}(s):
-\dfrac{\pi}{3}\le t-s,\ s\le \dfrac{\pi}{3}\}.\end{array}$$

Since $1= \max\{x:x\in P\cap\IR\}=\max\{{\rm Re\ } z:z\in P\}$, $1$ is the only possible star
center, if any, of $P$. Clearly, $0.95^2 i=\(0.95e^{i\frac{\pi}{4}}\)^2\in P$. From the following plot one can see that $t+(1-t)0.95^2 i\not\in P$ for $0<t<\dfrac{1}{3}$.
\fi

\begin{center}
\begin{tikzpicture}[line cap=round,line join=round,x=1.0cm,y=1.0cm,scale=2.5]
\draw[->,color=black] (-1,0.) -- (1.25,0.);
\draw[->,color=black] (0.,-1) -- (0.,1.25);
\draw [samples=50,variable=\y,domain=-0.8660254037844386: 0.8660254037844386] plot ({1/4-\y*\y},\y);
\draw (-0.5,-0.8660254037844386) -- (-0.2458780928473946,-0.9176295349746149) ;
\draw(-0.2458780928473946,-0.9176295349746149)--(0,-0.9025);
\draw (0,-0.9025) --(0.3223115666172949,-0.49913946280361987);
\draw(0.3223115666172949,-0.49913946280361987)--(0.9176295349746149,-0.2458780928473947) ;
\draw(0.9176295349746149,-0.2458780928473947) --(1,0);
\draw (1,0) --(0.9176295349746149,0.2458780928473947);
\draw (0.9176295349746149,0.2458780928473947)--(0.3223115666172949,0.49913946280361987);
\draw (0.3223115666172949,0.49913946280361987)--(0,0.9025);
\draw (0,0.9025) -- 
 (-0.2458780928473946,0.9176295349746149);
 \draw  (-0.2458780928473946,0.9176295349746149) --  
(-0.5,0.8660254037844386);
\draw (-0.5,-0.8660254037844386) -- (1,0);
\draw (-0.5,0.8660254037844386) -- (1,0);
\draw[dashed] (0,0.9025) -- (1,0);
\fill[fill opacity=0.2,samples=50,variable=\y,domain=0.8660254037844386:-0.8660254037844386] 
(-0.5,-0.8660254037844386)  -- (-0.2458780928473946,-0.9176295349746149)--(0,-0.9025) --(0.3223115666172949,-0.49913946280361987)--(0.9176295349746149,-0.2458780928473947) --(1,0) --(0.9176295349746149,0.2458780928473947)--(0.3223115666172949,0.49913946280361987)--
(0,0.9025) -- 
 (-0.2458780928473946,0.9176295349746149) --  
(-0.5,0.8660254037844386) --
plot ({1/4-\y*\y},\y) --  cycle;
\begin{scriptsize}
\fill (0.,0.9025) circle (1pt);
\draw (0.,0.9025) node[anchor=south west] {$\alpha_2^2$};
\fill (0.,-0.9025) circle (1pt);
\draw (0.,-0.9025) node[anchor=north west] {$\overline{\alpha_2}^2$};

\fill (0.9176295349746149,-0.2458780928473947) circle (1pt);
\draw (0.9176295349746149,-0.2458780928473947) node[right] {$\alpha_2\overline{\alpha_1}$};

\fill (0.9176295349746149,0.2458780928473947) circle (1pt);
\draw (0.9176295349746149,0.2458780928473947) node[right] {$\alpha_1\overline{\alpha_2}$};
\fill (-0.5,-0.8660254037844386) circle (1pt);
\draw (-0.5,-0.8660254037844386) node[below] {$\overline{\alpha_1}^2$};
\fill (1.,0.) circle (1pt);
\draw (1.1,0) node[above] {$1$};
\fill  (-0.2458780928473946,0.9176295349746149)  circle (1pt);
\draw  (-0.2458780928473946,0.9176295349746149)  node[above] {$\alpha_1\alpha_2$};
\fill  (-0.2458780928473946,-0.9176295349746149)  circle (1pt);
\draw  (-0.2458780928473946,-0.9176295349746149)  node[below] {$\overline{\alpha_1\alpha_2}$};

\fill  (-0.5,0.8660254037844386) circle (1pt);
\draw (-0.5,0.8660254037844386) node[above] {$\alpha_1^2$};
\end{scriptsize}
\end{tikzpicture}\medskip\\
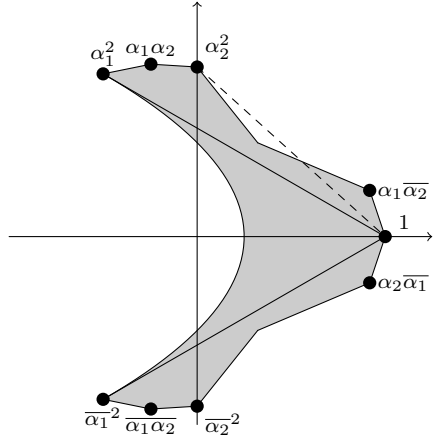
\textbf{Figure 4.} The set $P=K_1K_1$ in Example \ref{3.2} does not contain the segment $K(1,\alpha_2^2)$. 
\end{center}

\subsection{A necessary and sufficient condition}

In the following result, we establish a necessary and sufficient condition for the product of
two polygons to be a star-shaped set.

\begin{theorem} \label{3.3} Let
 $K_1=K(a_1,\ldots,a_n)$ and $K_2=K(b_1,\ldots,b_m)$. Then $K_1K_2$ is star-shaped if and only if there is $p \in K_1K_2$
such that $K(p,a_ib_j) \subseteq K_1K_2$ for all $1 \le i \le n$ and $1 \le j \le m$.
\end{theorem}

\it  Proof. \rm
Assume that $K_1=K(\a_1,\ldots,\a_n)$ and $K_2=K(\b_1,\ldots,\b_m)$.
From Proposition \ref{1.0} (a), we only need to prove that given any
$1\le i_1,i_2\le n$ and $1\le j_1,j_2\le m$, $K(p,q)\subseteq K_1K_2$ for all
$q\in K(\a_{i_1},\a_{i_2})K(\b_{j_1},\b_{j_2})$. Without loss of generality, we may assume that for $r=1,\ 2$, $i_r=j_r=r$, $\a_r=1+ia_r$ and $ \b_r=1+ib_r$ satisfy one of the conditions (a), (b) or (c) in Theorem \ref{thm2.4}.

Since   $K(p,\a_r\b_{t})\subseteq K_1K_2$ for $r,t=1,2$,
by the fact that $K_1K_2$ is
 simply connected, we see that
$$\mathbf{K}=K(p,\a_{1}\b_{1},\a_{1}\b_{2})\cup K(p,\a_{2}\b_{1},\a_{2}\b_{2})\cup
K(p,\a_{1}\b_{1},\a_{2}\b_{1})\cup K(p,\a_{1}\b_{2},\a_{2}\b_{2})\subseteq K_1K_2. $$

If $K(\a_{1},\a_{2})K(\b_{1},\b_{2})$ is convex, then $K(p,q)\subseteq \mathbf{K}$ for all $q\in K(\a_{1},\a_{2})K(\b_{1},\b_{2})$.

If $K(\a_{1},\a_{2})K(\b_{1},\b_{2})$ is not convex, then $a_1,\ a_2, b_1 $ and $ b_2$ satisfy conditions (b) or (c) in Theorem \ref{thm2.4}. Let $[a_1,a_2]\cap [b_1,b_2]=[c_1,c_2]$, $C=K(1+c_1i,1+c_2i)$, $B=K(1+b_1i,1+b_2i)\setminus C$ and $A=K(1+a_1i,1+a_2i)
\setminus C$. Since
$ K_1K_2= (A C)\cup (A B) \cup
(C  C) \cup (C  B) $,
and previous argument shows that $K(p,q)\subseteq K_1K_2$ for all $q\in (A  C)\cup (A B) \cup
 (C  B)$, it remains to show that  $K(p,q)\subseteq K_1K_2$ for all $q\in \partial (C  C)$. Let
$$\mathbf{V}= (1+c_1i)K(1+c_1i,1+c_2i)\cup (1+c_2i)K(1+c_1i,1+c_2i)\mbox{ and }\mathbf{U}=\{(1 + si)^2:s\in (c_1,c_2)\}\,.$$
 Note that $\partial (C  C)=\mathbf{V}\cup \mathbf{U} $ and $\mathbf{V}\subseteq K(\a_{1}\b_{1},\a_{1}\b_{2})\cup K(\a_{2}\b_{1},\a_{2}\b_{2})\cup
K(\a_{1}\b_{1},\a_{2}\b_{1})\cup K(\a_{1}\b_{2},\a_{2}\b_{2})$. So it remains to show that $K(p,q)\subseteq K_1K_2$ for all $q\in \mathbf{E}^o =\{(1 + si)^2:s\in (c_1,c_2)\}$.

Suppose $q\in \mathbf{E}^o$.
Let $\mathbf{L}$ be the tangent line to $\mathbf{E}^o$ at $q$   and $\mathbf{H}$  the open half plane determined by  $\mathbf{L}$ and contains $0$.  
\begin{center}
\begin{tikzpicture}[line cap=round,line join=round,x=1.0cm,y=1.0cm,scale=2]
\draw[->,color=black] (-1,0.) -- (1.2,0.);
\draw[->,color=black] (0.,-1) -- (0.,1.2);
\clip (-1,-1) rectangle (1.2,1.2);
\draw [samples=50,variable=\y,domain=-0.8660254037844386:0.8660254037844386] plot ({1/4-\y*\y},\y);
\draw (1.,0.)-- (-0.5,0.8660254037844386);
\draw (-0.5,-0.87)-- (1.,0.);
\draw[line width=1pt] (1,-1.5) -- (-0.5,1.73);
\fill[pattern=north west lines] (1,-1.5)-- ( -1,-1.5) -- (-1,1.5) -- (-0.5,1.73);
\fill[fill opacity=0.2,samples=50,variable=\y,domain=-0.8660254037844386:0.8660254037844386] plot ({1/4-\y*\y},\y) -- (1,0) -- cycle;
\begin{scriptsize}
\fill(0.19626790514376105,0.2318018439448637) circle (1pt);
\draw (0.19626790514376105,0.2318018439448637) node[right] {$q$};
\draw (-0.1,1.1) node {$\mathbf{L}$};
\draw (0.5,0.1) node {$CC$};
\end{scriptsize}
\fill[color=white] (-0.7,0.7) circle (3pt);
\draw (-0.7,0.7) node {$\mathbf{H}$};
\end{tikzpicture}\\
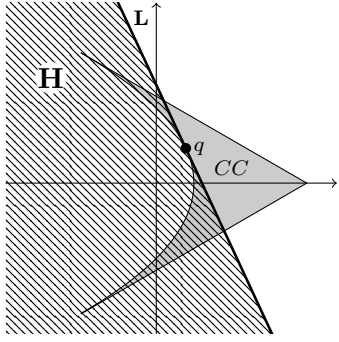
\textbf{Figure 5}
\end{center}

Consider the following three cases:

{\bf Case 1} If $p\in \mathbf{H}$, then there exists $t>1$ such that $s=p+t(q-p)\in \mathbf{V}$. Therefore, $  K(p,q)\subseteq K(p,s)\subseteq K_1K_2$.

{\bf Case 2} If $p\in \(\IC\setminus \mathbf{H}\)\cap (C  C)$, then $ K(p,q)\subseteq (C  C)\subseteq K_1K_2$ because $\(\IC\setminus \mathbf{H}\)\cap (C  C)$ is a triangular region containing $q$.

{\bf Case 3} If $p\in\IC\setminus  \(\mathbf{H}\cup (C  C)\)$, then there exists $0<t<1$ such that $s=p+t(q-p)\in \mathbf{V}$. Therefore, $  K(p,q)=K(p,s)\cup K(s,q)\subseteq  K_1K_2$.
\qed

\iffalse

\it  Proof. \rm
Assume that $K_1=K(a_1,\ldots,a_n)$ and $K_2=K(b_1,\ldots,b_m)$.
From Proposition \ref{1.0} (a), we only need to prove that given any
$1\le i_1,i_2\le n$ and $1\le j_1,j_2\le m$, $K(p,q)\subseteq K_1K_2$ for all
$q\in K(a_{i_1},a_{i_2})K(b_{j_1},b_{j_2})$.
Let $W$ be a subset of $\IC$  such that
$$\partial W=K(a_{i_1}b_{j_1},a_{i_1}b_{j_2})\cup K(a_{i_2}b_{j_1},a_{i_2}b_{j_2})
\cup K(a_{i_1}b_{j_1},a_{i_2}b_{j_1})\cup K(a_{i_1}b_{j_2},a_{i_2}b_{j_2}).$$
Since $\partial W\subseteq K_1K_2$ and $K(p,a_{i_r}b_{j_t})\subseteq K_1K_2$ for $r,t=1,2$,
by the fact that $K_1K_2$ is
 simply connected, we see that
$$K(p,a_{i_1}b_{j_1},a_{i_1}b_{j_2}),K(p,a_{i_2}b_{j_1},a_{i_2}b_{j_2}),
K(p,a_{i_1}b_{j_1},a_{i_2}b_{j_1}),K(p,a_{i_1}b_{j_2},a_{i_2}b_{j_2})\subseteq K_1K_2. $$
This implies that
\begin{equation}\label{eq-kz}
K(p,q)\subseteq K_1K_2\mbox{ for all }q\in W.
\end{equation}
On the other hand, $K(a_{i_1},a_{i_2})K(b_{j_1},b_{j_2})$ can
be written as
$$\left\{
      \begin{array}{cl}
        W & \mbox{ if }K(a_{i_1},a_{i_2})K(b_{j_1},b_{j_2})\mbox{ is convex} , \\
       W\cup E & \mbox{ otherwise, }
      \end{array} \right.$$
where $E=K(a_{i_1},a_{i_2})K(b_{j_1},b_{j_2})$;
see Theorem \ref{thm2.4} and its proof.
By (\ref{eq-kz}), we get the conclusion. \qed
\fi

We have the following consequence of Theorem \ref{3.3}.

\begin{corollary} \label{3.4}
Let $K_1$ be a triangle set with $K_1=\overline{K}_1$. Then $K_1=K(r,a,\overline{a})$
for some $r\in\IR $ and $a\in\IC$. The product set
$P=K_1K_1$ is a star-shaped set with $|a|^2$ as a star center.
\end{corollary}

\it Proof. \rm By Theorem \ref{3.3}, it suffices to show that $K(|a|^2,q)\in P$ for all $q\in \{r^2,ra,r\overline{a}, a^2,\overline{a}^2\}$.

\begin{enumerate}
\item
For $0\le t\le 1$, let $f(t)=(tr+(1-t)a)(tr+(1-t)\overline{a})\in P$. Since $f(0)=|a|^2$ and $f(1)=r^2$, we have $K(|a|^2,r^2)\in P$.

\item  $K(|a|^2,ra)=aK(\overline{a},r)\subseteq P$.

\item  $K(|a|^2,r(\overline{a})= \overline{a}K(a,r)\subseteq P$.

\item  $K(|a|^2,a^2)=aK(\overline{a},a)\subseteq P$.

\item  $K(|a|^2,\overline{a}^2)=\overline{a}K(a,\overline{a})\subseteq P$.
\end{enumerate}
\vskip -.35in  \qed

Suppose $A \in M_n$ is a real matrix. Then $W(A)$ is symmetric about the real axis.
By Corollary \ref{3.4}, if $A \in M_3$ is a real normal matrix, then
$W(A)W(A)$ is star-shaped. In fact,
if $A$ is Hermitian, then $W(A)W(A)$ is convex;
otherwise, $|a|^2$ is a star center, where $a, \bar a$ are the complex eigenvalues
of $A$.

\section{A line and a convex set}

In this section, we consider the product of a line segment and a convex set.
In the context of numerical range, we consider $W(A)W(B)$,
where  $A$ is a normal matrix with
collinear eigenvalues, and $B$ is a general matrix.

\begin{theorem}\label{star}
Let $K_1=K(\alpha,\beta)$ for some $\alpha, \beta\in \mathbb{C}$ and $K_2$ be a compact convex sets in $\mathbb{C}$. Then $K_1K_2$ is star-shaped.
\end{theorem}

We begin with the following easy cases.

\begin{proposition} \label{4.2}
Suppose that $K_1 = K(\a,\b)$ is a line segment and that $K_2$ is a (not necessarily compact)
convex set.
\begin{itemize}
\item[{\rm (1)}] If $0 \in K_1 \cup K_2$, then $K_1K_2$ is star-shaped with 0 as a star center.
\item
[{\rm (2)}] If there is a nonzero $\xi_1 \in \IC$ such that
$\xi_1 K_1 \subseteq (0, \infty)$,  then $K_1K_2$ is convex.
\item[
{\rm (3)}]
If there is a nonzero $\xi_1\in \IC$ such that $\xi_1 K_1\subseteq K_2$, then
$K_1K_2$ is star-shaped with $\xi_1 \a\b$ as a star center.
\end{itemize}
\end{proposition}

\it Proof. \rm  (1) It follows from Proposition \ref{1.0} (b).

(2) We may assume that $\xi_1 = 1$. Then
$K_1K_2 = \cup_{\a \le t \le \b} tK_2$ is convex.

(3) Assume $\xi_1=1$.
For every $p\in K_1$ and  $q\in K_2$, we will show that
$$K(\a\b,pq)\subseteq K(\a,\b)K(\a,\b,q)
\subseteq K_1K_2.$$
To this end, note that
$$K(\a\b,\a^2)=\a K(\a,\b),\ K(\a\b,\b^2)=\b K(\a,\b),\ K(\a\b,\a q)=\a K(\b,q),\ K(\a\b,\b q)=\b K(\a,q).$$
So, we have $K(\a\b, v) \in K(\a,\b)K(\a,\b,q)$
for any $v \in \{\a^2,\a\b, \a q, \b^2, \b q\}$,
which is the set of the product of vertexes of $K(\a,\b)$ and $K(\a,\b,q)$.
By Theorem \ref{3.3},
$K(\a,\b)K(\a,\b,q)$ is star-shaped with $\a\b$ as a star center.
Thus,
$K(\a\b, pq) \subseteq K(\a,\b)K(\a,\b,q) \subseteq K_1K_2.$

If $\xi_1 \ne 1$, then $(\xi_1\a)(\xi_1\b)$ is a star center of $(\xi_1 K_1)K_2 = \xi_1 K_1 K_2$
by the above argument.
Thus, $\xi_1(\a\b)$ is a star center of $K_1 K_2$.
\qed
\bigskip

From now on, we will focus on convex sets $K_1$ and $K_2$ that
 do not satisfy the
hypotheses in Proposition \ref{4.2} (1) -- (3).
In particular,   we may find $\xi_1$ and $\xi_2$ so that
$\xi_1 K_1 = K(\hat a, \hat b)$  and $\xi_2 K_2$ is a compact convex set
containing $\hat c, \hat d$ and lying in the cone
$$\cC = \{ t_1\hat c + t_2 \hat d:  t_1, t_2 \ge 0\},$$
where
$\hat a = 1+ia, \hat b = 1+ib, \hat c = 1+ic, \hat d = 1 + id$
with $a \le b, c \le d$. There could be five different configurations of the two
sets $\xi_1 K_1$ and $\xi_2 K_2$
as illustrated in Figure 6. (Here, we assume that Proposition \ref{4.2} (3) does not
hold so that we do not have the case $c\leq a< b \leq d$.)
If $K_1, K_2$ are put in these  ``canonical'' positions, we can describe the star centers
of $K_1K_2$ in the next theorem.
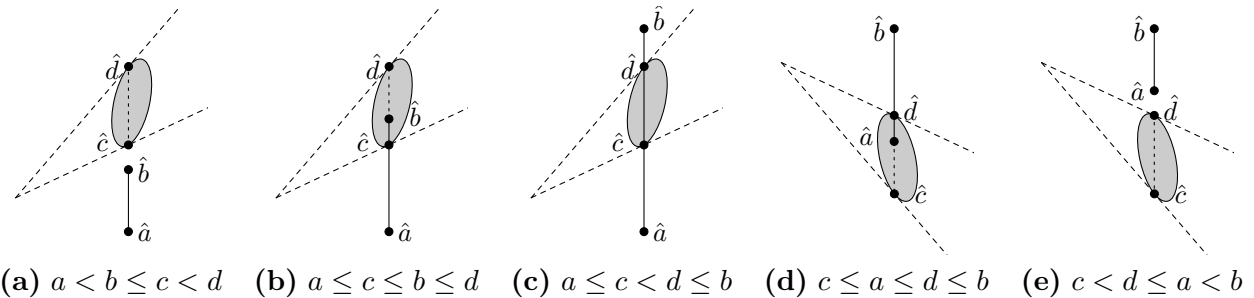
\begin{center}
\begin{tabular}{ccccc}
\begin{tikzpicture}[line cap=round,line join=round,x=1.0cm,y=1.0cm,scale=1.5]
\clip (-0.1,-0.5) rectangle (1.7,1.7);
\draw (1.,-0.3)-- (1.,0.25);
\draw[dash pattern=on 1pt off 2pt] (1.,0.46)-- (1.,1.165);
\draw [fill=black,fill opacity=0.2, rotate around={-102.9636070440955:(1.029059082343564,0.8427717986868933)}] (1.029059082343564,0.8427717986868933) ellipse (0.3970584957517196cm and 0.1569413367990951cm);
\draw [dash pattern=on 2pt off 2pt,domain=0.0:6.996459965006566] plot(\x,{(-0.--1.1647763366353459*\x)/1.0013392700642239});
\draw [dash pattern=on 2pt off 2pt,domain=0.0:6.996459965006566] plot(\x,{(-0.--0.46883474188652396*\x)/1.});
\begin{small}
\draw [fill=black] (1.,-0.3) circle (1pt);
\draw[color=black] (1,-0.3) node[right] {$\hat a$};
\draw [fill=black] (1.,0.25) circle (1pt);
\draw[color=black] (1,0.25) node[right] {$\hat b$};
\draw [fill=black] (1.0013392700642239,1.1647763366353459) circle (1pt);
\draw[color=black] (1,1.165) node[left] {$\hat d$};
\draw [fill=black] (1,0.46883474188652396) circle (1pt);
\draw[color=black] (0.9,0.5) node[left] {$\hat c$};
\end{small}
\end{tikzpicture}
\quad &
\begin{tikzpicture}[line cap=round,line join=round,x=1.0cm,y=1.0cm,scale=1.5]
\clip (-0.1,-0.5) rectangle (1.7,1.7);
\draw (1.,-0.3)-- (1.,0.7);
\draw[dash pattern=on 1pt off 2pt] (1.,0.7)-- (1.,1.165);
\draw [fill=black,fill opacity=0.2, rotate around={-102.9636070440955:(1.029059082343564,0.8427717986868933)}] (1.029059082343564,0.8427717986868933) ellipse (0.3970584957517196cm and 0.1569413367990951cm);
\draw [dash pattern=on 2pt off 2pt,domain=0.0:6.996459965006566] plot(\x,{(-0.--1.1647763366353459*\x)/1.0013392700642239});
\draw [dash pattern=on 2pt off 2pt,domain=0.0:6.996459965006566] plot(\x,{(-0.--0.46883474188652396*\x)/1.});
\begin{small}
\draw [fill=black] (1.,-0.3) circle (1pt);
\draw[color=black] (1,-0.3) node[right] {$\hat a$};
\draw [fill=black] (1.,0.7) circle (1pt);
\draw[color=black] (1.1,0.75) node[right] {$\hat b$};
\draw [fill=black] (1.0013392700642239,1.1647763366353459) circle (1pt);
\draw[color=black] (1,1.1647763366353459) node[left] {$\hat d$};
\draw [fill=black] (1,0.46883474188652396) circle (1pt);
\draw[color=black] (0.9,0.5) node[left] {$\hat c$};
\end{small}
\end{tikzpicture}\quad &
\begin{tikzpicture}[line cap=round,line join=round,x=1.0cm,y=1.0cm,scale=1.5]
\clip (-0.1,-0.5) rectangle (1.7,1.7);
\draw (1.,-0.3)-- (1.,1.5);
\draw [fill=black,fill opacity=0.2, rotate around={-102.9636070440955:(1.029059082343564,0.8427717986868933)}] (1.029059082343564,0.8427717986868933) ellipse (0.3970584957517196cm and 0.1569413367990951cm);
\draw [dash pattern=on 2pt off 2pt,domain=0.0:6.996459965006566] plot(\x,{(-0.--1.1647763366353459*\x)/1.0013392700642239});
\draw [dash pattern=on 2pt off 2pt,domain=0.0:6.996459965006566] plot(\x,{(-0.--0.46883474188652396*\x)/1.});
\begin{small}
\draw [fill=black] (1.,-0.3) circle (1pt);
\draw[color=black] (1,-0.3) node[right] {$\hat a$};
\draw [fill=black] (1.,1.5) circle (1pt);
\draw[color=black] (1,1.6) node[right] {$\hat b$};
\draw [fill=black] (1.0013392700642239,1.1647763366353459) circle (1pt);
\draw[color=black] (1,1.1647763366353459) node[left] {$\hat d$};
\draw [fill=black] (1,0.46883474188652396) circle (1pt);
\draw[color=black] (0.9,0.5) node[left] {$\hat c$};
\end{small}
\end{tikzpicture}\quad &
\begin{tikzpicture}[line cap=round,line join=round,x=1.0cm,y=-1.0cm,scale=1.5]
\clip (-0.1,-0.5) rectangle (1.7,1.7);
\draw (1.,-0.3)-- (1.,0.7);
\draw[dash pattern=on 1pt off 2pt] (1.,0.7)-- (1.,1.165);
\draw [fill=black,fill opacity=0.2, rotate around={102.9636070440955:(1.029059082343564,0.8427717986868933)}] (1.029059082343564,0.8427717986868933) ellipse (-0.3970584957517196cm and -0.1569413367990951cm);
\draw [dash pattern=on 2pt off 2pt,domain=0.0:6.996459965006566] plot(\x,{(-0.--1.1647763366353459*\x)/1.0013392700642239});
\draw [dash pattern=on 2pt off 2pt,domain=0.0:6.996459965006566] plot(\x,{(-0.--0.46883474188652396*\x)/1.});
\begin{small}
\draw [fill=black] (1.,-0.3) circle (1pt);
\draw[color=black] (1,-0.3) node[left] {$\hat b$};
\draw [fill=black] (1.,0.7) circle (1pt);
\draw[color=black] (0.9,0.65) node[left] {$\hat a$};
\draw [fill=black] (1.0013392700642239,1.1647763366353459) circle (1pt);
\draw[color=black] (1.1,1.1647763366353459) node[right] {$\hat c$};
\draw [fill=black] (1,0.46883474188652396) circle (1pt);
\draw[color=black] (1,0.4) node[right] {$\hat d$};
\end{small}
\end{tikzpicture}
 \quad &
\begin{tikzpicture}[line cap=round,line join=round,x=1.0cm,y=-1.0cm,scale=1.5]
\clip (-0.1,-0.5) rectangle (1.7,1.7);
\draw (1.,-0.3)-- (1.,0.25);
\draw[dash pattern=on 1pt off 2pt] (1.,0.46)-- (1.,1.165);
\draw [fill=black,fill opacity=0.2, rotate around={102.9636070440955:(1.029059082343564,0.8427717986868933)}] (1.029059082343564,0.8427717986868933) ellipse (-0.3970584957517196cm and -0.1569413367990951cm);
\draw [dash pattern=on 2pt off 2pt,domain=0.0:6.996459965006566] plot(\x,{(-0.--1.1647763366353459*\x)/1.0013392700642239});
\draw [dash pattern=on 2pt off 2pt,domain=0.0:6.996459965006566] plot(\x,{(-0.--0.46883474188652396*\x)/1.});
\begin{small}
\draw [fill=black] (1.,-0.3) circle (1pt);
\draw[color=black] (1,-0.3) node[left] {$\hat b$};
\draw [fill=black] (1.,0.25) circle (1pt);
\draw[color=black] (1,0.25) node[left] {$\hat a$};
\draw [fill=black] (1.0013392700642239,1.1647763366353459) circle (1pt);
\draw[color=black] (1.1,1.1647763366353459) node[right] {$\hat c$};
\draw [fill=black] (1,0.46883474188652396) circle (1pt);
\draw[color=black] (1,0.4) node[right] {$\hat d$};
\end{small}
\end{tikzpicture} \\
\textbf{(a)}  $a<b\leq  c<d$  & \textbf{(b)} $a\leq c\leq b\leq  d$ & \textbf{(c)} $a\leq c < d\leq b$ & \textbf{(d)} $c\leq a \leq d\leq b$  & \textbf{(e)} $c< d \leq a< b$
\end{tabular}\medskip\\
\textbf{Figure 6.} The above figures illustrate the canonical representations of a line segment $K_1=K(a,b)$ and a convex set $K_2$ described in Theorem  \ref{star-2}.
\end{center}

\begin{theorem} \label{star-2}
Let $\hat a = 1+ia, \hat b = 1+ib, \hat c = 1+ic, \hat d = 1 + id$
with $a \le b, c \le d$.
Suppose $K_1 = K(\hat a, \hat b)$  and $K_2$ be a compact convex set
containing $\hat c, \hat d$ and lying in the cone
$$\cC = \{ t_1\hat c + t_2 \hat d:  t_1, t_2 \ge 0\}$$
such that the hypotheses of Proposition \ref{4.2} {\rm (1) -- (3)} do not hold.
Then $K_1K_2$ is star-shaped and one of the following holds.

{\rm (a)} If $a \le b \le c \le d$, then
 $\hat b \hat c$ is a star center.

{\rm (b)} If $a \le c \le b \le d$, then  $\hat b \hat c$ is a star center.

{\rm (c)} If $a \le c \le d \le b$, then $\hat c\hat d$ is a star center.

{\rm (d)} If $c \le a \le d \le b$, then  $\hat a \hat d$ is a star center.

{\rm (e)} If $c \le d \le a \le b$, then
$\hat a \hat d $ is a star center.

\end{theorem}

We need some lemmas to prove Theorem \ref{star-2}.

\begin{lemma} \label{cp} Suppose $C=1+i\tan\tc$, $D=1+i\tan\td$  and $P=re^{i\tp}$ with $r>0$, $-\dfrac{\pi}{2}<\tc< \tp< \td<\dfrac{\pi}{2}$.
Let
$$\dfrac{-i(P-C)}{|P-C|}=e^{i\t_1}\ \mbox{ and }\ \dfrac{i(P-D)}{|P-D|}=e^{i\t_2}\ \mbox{ with }\ -\dfrac{\pi}{2}<\t_1,\ \t_2< \dfrac{\pi}{2}\,. $$ Then there exists $\xi_1,\ \xi_2$ such that $ \xi_1C =1+i\tan (\tc-\t_1)$ and $ \xi_1P =1+i\tan(\tp-\t_1)$,  $ \xi_2D =1+i\tan (\td-\t_2)$ and $ \xi_2P =1+i\tan(\tp -\t_2)$.

Consequently, we have

\begin{enumerate}
\item If $\Re(P)\le 1$, then $\t_2\le 0\le \t_1$ and $\tc-\t_1\le \tp-\t_1\le \tp\le \tp-\t_2\le   \td\le \td-\t_2$.
\item If $\Re(P)\ge 1$, then $\t_1\le 0\le \t_2$ and $\tc\le \tc-\t_1\le \tp-\t_1$ and $   \tp-\t_2\le   \td-\t_2\le \td$.
    \end{enumerate}
\end{lemma}

\it Proof. \rm First consider $C$ and $P$.
Then $\t_1$ is the angle from  $\overrightarrow{CD}$ to   $\overrightarrow{CP}$. Then the result follows from simple geometry.\vskip.2in

\begin{center}
\begin{tikzpicture}[line cap=round,line join=round,x=1.0cm,y=0.7cm,scale=2]
\draw[->] (-0.2,0.) -- (2,0.);
\draw[->] (0.,-0.2) -- (0.,3.5);

\draw [shift={(0.,0.)},fill=black,fill opacity=0.1] (0,0) -- (0.:0.443110057396134) arc (0.:18.434948822922006:0.443110057396134) -- cycle;
\draw[fill=black,fill opacity=0.1] (1.2203670316059365,0.5388989051821906) -- (1.1014681264237458,0.49926593678812714) -- (1.1411010948178093,0.3803670316059364) -- (1.26,0.42) -- cycle; 
\draw [shift={(1.,1.2)},fill=black,fill opacity=0.1] (0,0) -- (90.:0.443110057396134) arc (90.:108.43494882292201:0.443110057396134) -- cycle;
\draw [shift={(1.,3.3)},fill=black,fill opacity=0.1] (0,0) -- (-104.0362434679265:0.443110057396134) arc (-104.0362434679265:-90.:0.443110057396134) -- cycle;
\draw [shift={(0.,0.)},fill=black,fill opacity=0.1] (0,0) -- (0.:0.2658660344376804) arc (0.:71.56505117707799:0.2658660344376804) -- cycle;
\draw [shift={(0.,0.)},fill=black,fill opacity=0.1] (0,0) -- (0.:0.35448804591690725) arc (0.:50.19442890773481:0.35448804591690725) -- cycle;
\draw[dashed] (0,0) -- (1,3.3);
\draw (1.,3.3)-- (1.,1.2);
\draw (0.,0.)-- (0.7,2.1);
\draw (0.7,2.1)-- (1.,1.2);
\draw (0.,0.)-- (1.,1.2);
\draw (1.4,0.)-- (0.7,2.1);
\draw (1.26,0.42)-- (0.,0.);
\draw (0.7,2.1)-- (1.,3.3);
\begin{scriptsize}
\draw [fill=black] (1.,3.3) circle (1pt);
\draw[color=black] (1.1453148428626505,3.2405634850878355) node {$D$};
\draw [fill=black] (1.,1.2) circle (1pt);
\draw[color=black] (1.1541770440105732,1.1845328187697794) node {$C$};
\draw [fill=black] (0.,0.) circle (1pt);
\draw[color=black] (-0.05994451325483416,-0.1270729511227737) node {$0$};
\draw [fill=black] (0.7,2.1) circle (1pt);
\draw[color=black] (0.8085511992415886,2.168237146189194) node {$P$};
\draw [fill=black] (1.26,0.42) circle (1pt);
\draw[color=black] (1.3757320727086402,0.5464543361193482) node {$R$};
\draw [fill=black] (1.4,0.) circle (1pt);
\draw[color=black] (1.5175272910754032,0.12106868101906071) node {$X$};
\draw[color=black] (0.6,0.076) node {$\theta_1$};
\draw[color=black] (0.94,1.74) node {$\theta_1$};
\draw[color=black] (0.94,2.7) node {$\theta_2$};
\draw[color=black] (0.258,0.45) node {$\theta_P$};
\draw[color=black] (0.45,0.265) node {$\theta_C$};
\end{scriptsize}
\end{tikzpicture}\\
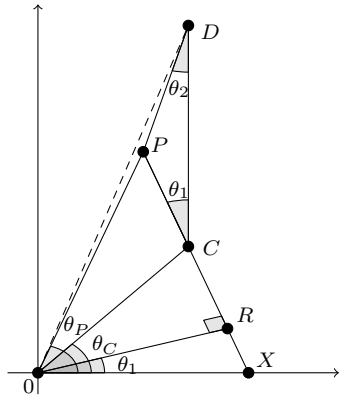
\textbf{Figure 7.}
\end{center}

\vskip.2in
On one also can calculate directly with $\xi_1=\dfrac{\cos\tc}{\cos(\tc-\t_1)}e^{-i\t_1}$.

For the second statement, apply the above result on $\overline{D}$ and $\overline{P}$, the complex conjugate of $D$ and $P$.
\qed

\begin{lemma}\label{lem0.1} Suppose $a \le c \le d$,
$p = t_1 (1+ic) + t_2(1+id)$ is nonzero for some $t_1, t_2 \ge 0$,
$K_1 = K(1+ia,1+id)$, and $K_2 = K(1+ic, 1+id, p)$.
Then $K_1K_2$ is star-shaped with $(1+ic)(1+id)$ as a star center.
\end{lemma}

\it Proof. \rm  Let $\hat a = 1+ia$, $\hat c = 1+ic$, $\hat d = 1 + id$.
By Theorem \ref{3.3}, it suffices to show that $K(\hat c \hat d,uv) \subseteq K_1K_2$
for each pair of elements
$(u,v)$ in $\{\hat{a},\hat{d}\} \times \{\hat{c}, \hat{d}, p\}.$ If $u=\hat{d}$,
then $K(\hat{c}\hat{d},\hat{d}v)=\hat{d}K(\hat{c},v)\subseteq K_1K_2$. Similarly,
if $u=\hat{c}$, then $K(\hat{c}\hat{d},\hat{c}v)=\hat{c}K(\hat{d},v)\subseteq K_1K_2$.
Thus, the only nontrivial case is when $(u,v) = (\hat{a}, p)$.

By continuity, we may assume that $t_1,\ t_2>0$. We consider two cases.

{\bf Case 1} Suppose $\Re(p) \le 1$.
Then by Lemma \ref{cp} and Theorem \ref{thm2.4}, $ K(\hat a, \hat c) K(p,\hat d)$ is convex. So

$$K(\hat c \hat d,\hat{a}p) \subseteq K(\hat a, \hat c) K(p,\hat d) \subseteq  K_1K_2\,.$$

\medskip
{\bf Case 2} Suppose  $\Re(p) > 1$. By Lemma \ref{cp}, there exists $\a_0$ such that $\a_0\hat c=1+c_1i$ and $\a_0p=1+p_1i$ such that $c_1>c$. By Theorem \ref{thm2.4}, if $p_1\ge d$, then  $\hat c\hat d$ is a star center of $K(\hat a,\hat d)K(\hat c, p)$. If $p_1< d$, then $K(\hat a\hat c,\hat d\hat c)$ intersects $K(\hat ap,\hat d p)$ and $\hat c \hat d$ lies inside the triangle with vertices $\hat ap,\ \hat d p,\ \hat a  \hat d$ (see Figure 8). Thus, $K(\hat{c}\hat{d},\hat{a}p)$ is in the interior of the region enclosed by $K(\hat{d}p,\hat{c}\hat{d})\cup K(\hat{c}\hat{d}, \hat{a}\hat{d}) \cup K(\hat{a}\hat{d},\hat{a}p)\cup K(\hat{a}p,\hat{c}\hat{a})\subseteq K_1K_2$. 
\medskip
\begin{center}
\begin{tikzpicture}[line cap=round,line join=round,x=1.5cm,y=0.7cm]
\draw (-1.25,3.)-- (1.25,1.875);
\draw[line width=1.5pt] (2.5,0.5)-- (1.25,-1.25);
\draw[line width=1.5pt] (1.25,-1.25)-- (1.25,1.875);
\draw[line width=1.5pt,dashed] (2.5,0.5)-- (1.25,1.875);
\draw (2.5,0.5)-- (-1.25,3.);
\draw (1.75,1)-- (2.5,0.5);
\draw (1.75,1.)-- (0.5,-1.5);
\draw (1.25,-1.25)-- (0.5,-1.5);
\draw (1.75,1.)-- (1.25,1.875);
\draw (0.5,-1.5)-- (2.5,0.5);
\draw [samples=50,domain=-1:3,variable=\y] plot ({1-(\y*\y)/4},\y);
\fill[fill opacity=0.2,samples=50,domain=-1:3,variable=\y] plot ({1-(\y*\y)/4},\y)--(1.25,1.875)-- (1.75,1) -- (2.5,0.5) --(1.25,-1.25) -- (0.5,-1.5) --  cycle;
\begin{scriptsize}
\fill (1.25,1.875) circle (1.5pt);
\draw (1.25,1.875) node[above] {$\hat{d}p$};
\fill (1.25,-1.25) circle (1.5pt);
\draw (1.25,-1.25) node[anchor=north west] {$\hat{a}p$};
\fill  (-1.25,3) circle (1.5pt);
\draw (-1.25,3) node[left] {$\hat{d}^2$};
\draw (2.5,0.5) node[right] {$\hat{a}\hat{d}$};
\fill (2.5,0.5) circle (1.5pt);
\fill (1.75,1) circle (1.5pt);
\draw (1.75,1) node[below] {$\hat{c}\hat{d}$};
\fill (0.5,-1.5) circle (1.5pt);
\draw (0.5,-1.5) node[below] {$\hat{c}\hat{a}$};
\fill (0.75,-1) circle (1.5pt);
\draw (0.75,-1) node[left] {$\hat{c}^2$};
\end{scriptsize}
\end{tikzpicture}\\
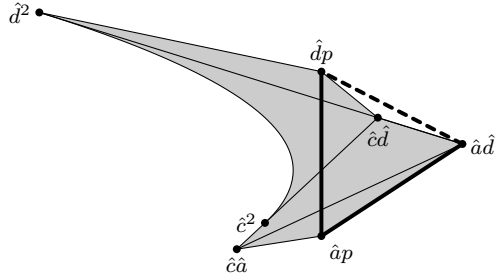
\textbf{Figure 8.}
\end{center}
\medskip
 In both cases, we have
$K(\hat c \hat d,\hat ap)\subseteq K_1K_2$.\qed

\begin{lemma}\label{lem0.2}
Suppose $a < b \le c < d$,
$p = t_1 (1+ic) + t_2(1+id)$ is nonzero for some $t_1, t_2 \ge 0$ and
$K_1 = K(1+ia,1+ib)$, and $K_2 = K(1+ic, 1+id, p)$. Assume also that there is no $\xi\in\IC$
such that $K_1\subseteq \xi K_2$. Then $K_1K_2$ is star-shaped and
$(1+bi)(1+ci)$ is a star center.
\end{lemma}

\it Proof. \rm  Let $\hat a = 1+ia$, $\hat c = 1+ic$, $\hat d = 1 + id$. Similar to the previous lemma, it is enough to show that $K(\hat{b}\hat{c},\hat{a}p)\subseteq K_1K_2$ for any $p=t_1\hat{c}+t_2\hat{d}$ such that $t_1,t_2\geq 0$.

Let $\xi\in\IC$
such that  $\xi K(\hat{c},p)$ is a vertical line segment with real part $1$.
If $\xi K(\hat{c},p)\not\subseteq K(\hat{a},\hat{b})$, then by Corollary \ref{2.6},
$\hat{b}\hat{c}$ is a star-center of $K_1K(\hat{c},p)$ and hence
$K(\hat{b}\hat{c},\hat{a}p)\subseteq K_1K_2$. Otherwise, we have
$\xi K(\hat{c},p)\subseteq K(\hat{a},\hat{b})$ and $K_1K(\hat{c},p)$ is as shown in Figure 9(c)
in the figure below. This will only happen if $\Re(p)<1$. Since
 $\hat{a}p=t_1(\hat{c}\hat{a})+t_2\hat{d}\hat{a}$ for some $t_1,t_2\geq 0$ such that $t_1+t_2<1$, then
$\hat{a}p\in K(0,\hat{c}\hat{a},\hat{d}\hat{a})$ and
$\hat{b}p\in K(0,\hat{c}\hat{b},\hat{d}\hat{b})$. Note also that $0$ and $p\hat{a}$
are separated by the line segment $K(\hat c\hat{b},\hat{c}\hat{a})$.
Hence, $p\hat{a}$ is in the quadrilateral $K_1K(\hat{c},\hat{d})$ and therefore
 $K(\hat{a}p,\hat{c}\hat{b})\subseteq K_1K_2$. This finishes the proof that $\hat{c}\hat{b}$ is a
 star center for $K_1K_2$. \qed
\begin{center}
\begin{tabular}{ccc}
\begin{tikzpicture}[line cap=round,line join=round,x=1.0cm,y=1.0cm,scale=1.5]
\fill[fill=black,fill opacity=0.2] (1.,1) -- (1.,0.5) -- (0.8,0.7) -- cycle;
\draw[dashed] (0,0)-- (1.5,1.5);
\draw[dashed] (0,0)-- (1.5,0.75);
\draw (1.,1)-- (1.,0.5);
\draw (1.,0.5)-- (0.8,0.7);
\draw (0.8,0.7)-- (1.,1);
\draw (1.,0.4)-- (1.,-1.);
\begin{scriptsize}
\fill  (1,0.4) circle (1pt);
\draw (1,0.4) node[right] {$\hat{b}$};
\fill  (1,-1) circle (1pt);
\draw (1,-1) node[right] {$\hat{a}$};
\fill  (1,1) circle (1pt);
\draw (1,1) node[left] {$\hat{d}$};
\fill  (1,0.5) circle (1pt);
\draw (0.95,0.6) node[right] {$\hat{c}$};
\fill (0.8,0.7) circle (1pt);
\draw (0.8,0.7) node[anchor= north east] {$p$};
\fill (0,0) circle (1pt);
\draw (0,-0.1) node[below] {$0$};
\end{scriptsize}
\end{tikzpicture} &
\begin{tikzpicture}[line cap=round,line join=round,x=1.0cm,y=1.0cm,scale=1.7]
\fill[fill=black,fill opacity=0.2] (0.6,1.4) -- (2,0) -- (1.5,-0.5) -- (0.8,0.9) -- cycle;
\draw (0.6,1.4) -- (2,0);
\draw (2,0) -- (1.5,-0.5);
\draw (1.5,-0.5) -- (0.8,0.9);
\draw (0.6,1.4) -- (0.8,0.9);
\draw[dashed] (0,0) -- (1.5,-0.5);
\draw[dashed] (0,0) -- (2,0);
\fill[pattern=north west lines] (1.25,0) -- (2,0) -- (2,0) -- (1.5,-0.5) -- cycle;
\begin{scriptsize}
\fill  (0,0) circle (1pt);
\draw (0,0) node[above] {$0$};
\fill[color=white]  (1.55,-0.2) circle (4pt);
\draw (1.55,-0.2) node {$p\hat{a}$};
\fill  (0.6,1.4) circle (1pt);
\draw (0.6,1.4) node[above] {$\hat{d}\hat{b}$};
\fill (2,0) circle (1pt);
\draw (2,0) node[right] {$\hat{d}\hat{a}$};
\fill  (1.5,-0.5) circle (1pt);
\draw (1.5,-0.5) node[left] {$\hat{c}\hat{a}$};
\fill  (0.8,0.9) circle (1pt);
\draw (0.8,0.9) node[left] {$\hat{c}\hat{b}$};
\end{scriptsize}
\end{tikzpicture} \qquad
& \qquad
\begin{tikzpicture}[line cap=round,line join=round,x=1.0cm,y=1.0cm,scale=2]
\draw [samples=50,domain=-0.33:-0.067,variable=\y] plot ({3/4-3*\y/2-3*\y*\y/4},{3*\y/2-3*\y*\y/4+3/4});
\fill[fill=black,fill opacity=0.2] (0.52,1.02) -- (1.5,-0.1) -- (1.5,-0.5) -- (0.8,0.9) -- cycle;
\draw (0.52,1.02) -- (0.8,0.9);
\draw (1.5,-0.1) -- (1.5,-0.5);
\draw (1.5,-0.5) -- (0.8,0.9);
\draw (0.52,1.02) -- (1.5,-0.1);
\begin{scriptsize}
\fill  (0.52,1.02) circle (1pt);
\draw (0.52,1.02) node[above] {$p\hat{b}$};
\fill (1.5,-0.1) circle (1pt);
\draw (1.5,-0.1) node[right] {$p\hat{a}$};
\fill  (1.5,-0.5) circle (1pt);
\draw (1.5,-0.5) node[left] {$\hat{c}\hat{a}$};
\fill  (0.8,0.9) circle (1pt);
\draw (0.8,0.9) node[right] {$\hat{c}\hat{b}$};
\end{scriptsize}
\end{tikzpicture}\\
\textbf{(a)} $K_1=K(\hat a,\hat{b})$ and $K(\hat{c},\hat{d},p)$ & \textbf{(b)} $K_1K(\hat{c},\hat{d})$ & \textbf{(c)} $K_1K(\hat{c},p)$
\end{tabular}\medskip\\
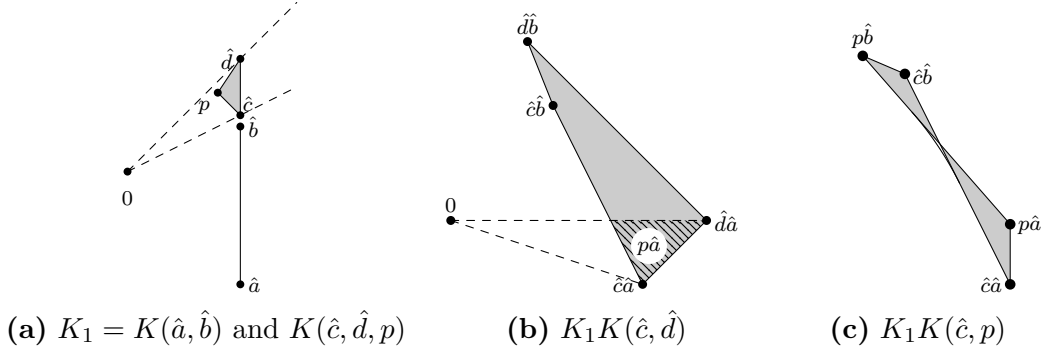
\textbf{Figure 9.}
\end{center}
\qed

\medskip
{\it Proof of Theorem \ref{star-2}:}
Note that (d) follows from (b) by considering $\bar{K}_1\bar{K}_2$. Similarly,
(e) follows from (a).   Thus, we only need to prove (a)-(c).

To prove that $s$ is a star center of $K_1K_2$, we show that for any
$p\in K_2$, $s$ is a star center of $K_1K(\hat{c},\hat{d},p)$. To accomplish this, it is
enough to show that $K(s,uv)\subseteq K_1K_2$ for all pairs
$(u,v)\in  \{\hat{b},\hat{a} \}\times \{\hat{c},\hat{d},p\}$ by Theorem \ref{3.3},
where $p=t_1\hat{c}+t_2\hat{d}$
for some $t_1,t_2\geq 0$.

For (a), the conclusion follows directly from Lemma \ref{lem0.2}.

To prove (c), the only nontrivial cases to consider are when $(u,v)=(\hat{a},p)$ or $(u,v)=(\hat{b},p)$. By Lemma \ref{lem0.1}, $K(\hat{c}\hat{d},\hat{a}p)\subseteq K(\hat{a},\hat{d})K(\hat{c},\hat{d},p)\subseteq K_1K_2$.
By Lemma \ref{lem0.1} again,
the product $\overline{K(\hat{b},\hat{c})}\ \overline{K(\hat{c},\hat{d},\hat{p})}$, has $\overline{\hat{c}\hat{d}}$ as a star center. Thus, $\hat{c}\hat{d}$ is a star center of $K(\hat{b},\hat{c})K(\hat{c},\hat{d},\hat{p})$ and thus $K(\hat{c}\hat{d},\hat{b}p)\subseteq K(\hat{b},\hat{c})K(\hat{c},\hat{d},\hat{p})\subseteq K_1K_2$.

To prove (b), it is enough to show that $K(\hat{c}\hat{b},\hat{a}p)\subseteq K_1K_2$ for all $p\in K_2$. We consider two cases,

\begin{enumerate}
\item[1.]Suppose  $p=t_1\hat{d}+t_2\hat{b}$ for some $t_1,t_2\geq 0$. Then by Lemma
\ref{lem0.2}, $\hat{b}\hat{c}$ is a star-center of $K(\hat{a},\hat{c})K(\hat{b},\hat{d},p)$. Thus $K(\hat{b}\hat{c},\hat{a}p)\subseteq K(\hat{a},\hat{c})K(\hat{b},\hat{d},p) \subseteq K_1K_2$.
\item[2.] Suppose  $p=t_1\hat{b}+t_2\hat{c}$ for some $t_1,t_2\geq 0$. Then by Lemma
\ref{lem0.1}, $\hat{b}\hat{c}$ is a star-center of $K(\hat{a},\hat{b})K(\hat{b},\hat{c},p)$. Thus $K(\hat{b}\hat{c},\hat{a}p)\subseteq K(\hat{a},\hat{b})K(\hat{b},\hat{c},p)\subseteq K_1K_2$.
\end{enumerate}
In both cases, $\hat{b}\hat{c}$ is a star-center for $K_1K_2$.
\qed

It is clear that Theorem \ref{star} follows from Proposition \ref{4.2}
and Theorem \ref{star-2}.

\section{A circular disk and a closed set}

It is known that the product of two circular disks is star-shaped \cite{FMR, FP, Mc, Karol}.
In this section, we will prove  some unexpected results that
if $K_1$ is a circular disk, then for many closed sets $K_2$, the product set
is star-shaped. We will use $D(\mu, R)$ to denote the closed disk with center $\mu\in\IC$ and
radius $R\ge 0$.

Note that if $0\in K_1$, then for every non-empty set $K_2$, $K_1K_2$ is star-shaped with $0$
as star center. Suppose $0\not\in K_1$, we can always scale $K_1$ so that it is a circular
disk centered at 1 with radius $r < 1$.

We have the following results showing that the product set of a circular disk
and another set would be star-shaped under some very general conditions.
We begin with the following observation.

\begin{lemma} \label{5.2}
Suppose   $r\in (0, 1]$
and  $b\in D(1,r)$. Then the product $D(1,r)\{b\}$ is a disk containing $1-r^2$.
\end{lemma}

\it Proof. \rm Let $b\in  D(1,r)$. Then $bD(1,r)=D(b,|b|r)$.
\begin{eqnarray*}
|b-(1-r^2)|^2&=&(b-(1-r^2))(\overline{b}-(1-r^2))\\
&=&|b|^2-(b+\overline{b})(1-r^2)+(1-r^2)^2\\
&=&|b|^2r^2-(1-r^2)(-|b|^2+(b+\overline{b})-(1-r^2))\\
&=&|b|^2r^2-(1-r^2)(r^2-(b-1)(\overline{b}-1))\\
&\le&|b|^2r^2\hskip .7in \mbox{because }|b-1|\le r\le 1.
\end{eqnarray*}
\vskip -.35in \qed

\medskip
From the above simple proposition, we get the following.

\begin{theorem} \label{5.1}  Suppose $K_1=D(\mu, R)$ does  not contain 0.
For every nonempty subset $S$ of $K_1$, the product set $K_1S$ is star shaped
with star center $\mu^2(1-r^2)$, where $r =|\mu^{-1}R|$.
\end{theorem}

In the numerical range context, for every circular disk $K_1$, there is
$A \in M_2$ such that $A - (\tr A)I/2$ is nilpotent and $W(A) = K_1$.
Moreover, $B \in M_n$ satisfies  $W(B) \subseteq W(A)$
if and only if $B$ admits a dilation of the form $I \otimes A$; see \cite{ando,CL}.
By Theorem \ref{5.1}, if $A \in M_2$
such that $(A-\tr A I)/2$ is nilpotent, then $W(A)W(B)$ is star-shaped
for any $B \in M_n$ satisfying $W(B) \subseteq W(A)$.

Next, we have the following.

\begin{theorem} \label{thmb}
Suppose $r\in (0,1]$ and $b \in \IC$ with $\Re(b) \ge 1$.
%$K_1$ is the close disk centered at $1$ with radius
%=1+Re^{i\theta}$, with $-\dfrac{\pi}2\le \theta\le \dfrac{\pi}2$.
Then the product $K(1,b)D(1,r)$ is  star-shaped with $1$ as star center.
\end{theorem}

\it Proof. \rm Suppose $b=Re^{i\theta}$
with $R\ge 0$ and $-\frac{\pi}{2}\le \theta\le \frac{\pi}{2}$.
Let $c\in K(1,b)$. Then $c=1+sRe^{i\theta}$ for  some $0\le s\le 1$. $cK_1=D(c,|c|r)$.
Therefore,
$K(1,b)D(1,r)= \cup\{ D(c,|c|r):c\in K(1,b)\}$. Let $z\in K(1,b)D(1,r)$. Then
$|z-(1+sRe^{i\theta})|\le  |1+sRe^{i\theta} |r$ for some  $0\le s\le 1$. Let $0\le t\le 1$.
We have
\begin{eqnarray*}
&&|tz+(1-t)-(1+tsRe^{i\theta})|^2\\
&=&|t(z-(1+sRe^{i\theta}))|^2\\
&\le&t^2|1+sRe^{i\theta} |^2r^2\\
&=&\((t+tsR\cos\theta)^2+(tsR\sin \theta)^2\)r^2\\
&=&\((1+tsR\cos\theta)^2+(tsR\sin \theta)^2-(1-t)(1+t+2tsR\cos\theta )\)r^2\\
&\le&\((1+tsR\cos\theta)^2+(tsR\sin \theta)^2\)r^2\\
&=&|1+tsRe^{i\theta}|^2r^2.
\end{eqnarray*}
Therefore, $tz+(1-t)\in D(1+tsRe^{i\theta},|1+tsRe^{i\theta}|r)\subseteq K(1,b)D(1,r)$.
\qed

\begin{theorem}  \label{thms} Suppose $S$ is a star-shaped subset of $\IC$ with star
center $s$ such that
$|s|\le |z|$ for  every $z\in S$. Then $D(a,r)S$ is star-shaped for every  circular disk
$D(a,r)$.  In particular, if $S$ is convex, then $D(a,r)S$ is star-shaped for every
circular disk $D(a,r)$.
\end{theorem}

\it Proof. \rm If either $S$ or  $D(a,r)$ contains $0$, the result holds. So we may assume that
$0\not\in S\cup D(a,r)$.

We may assume that $s=1$ and $D(a,r)=D(1,r) $ with $0\le r\le 1$. Then for every $z\in S$,
$z=1+Re^{i\theta}$ for some $-\dfrac{\pi}2\le \theta\le \dfrac{\pi}2$. By Theorem \ref{thmb},
the product $K(1,z)D(1,r)$ is star shaped with star center $1$. Hence, $SD(1,r)$ is also star
shaped with star center 1.
\qed

Apart from the nice results above, there are some limitations  about the
star-shapedness of the product set of
a circular disk and another set in $\IC$ as shown in the following.

\begin{example} \rm
Let $S=K(1,2e^{i\frac{11\pi}{12}})\cup K(1,2e^{-i\frac{11\pi}{12}})$. Then $S$ is
star-shaped with $1$ as star center. Let $D(1,\frac{1}{2})$ be the disk centered at
1 with radius $\frac{1}{2}$. Then the product set
$SD(1,\frac{1}{2}) $ is not simply connected.

\begin{center}
\begin{tikzpicture}[line cap=round,line join=round,x=1.0cm,y=1.0cm,scale=1.5]
\draw[<->,color=black] (-3.4,0.) -- (2,0.);
\draw[<->,color=black] (0.,-1.8) -- (0.,1.8);

%CIRCLES
\draw(1.,0.) circle (0.5cm);
\draw (-1.93185,0.51764) circle (1.cm);
\draw [samples=50,domain=-0.81:0.91,rotate around={-100.02530678696402:(0.03520175882521048,0.1689626876062166)},xshift=0.03520175882521048cm,yshift=0.1689626876062166cm] plot ({0.08319873700947768*(-1-(\x)^2)/(1-(\x)^2)},{0.1447370102073152*(-2)*(\x)/(1-(\x)^2)});
\draw [samples=50,domain=-0.91:0.81,rotate around={-100.02530678696402:(0.03520175882521048,0.1689626876062166)},xshift=0.03520175882521048cm,yshift=0.1689626876062166cm] plot ({0.08319873700947768*(1+(\x)^2)/(1-(\x)^2)},{0.1447370102073152*2*(\x)/(1-(\x)^2)});

\fill[fill opacity=0.3] (1,0) circle (0.5cm) --
plot[samples=50,domain=-0.81:0.91,rotate around={-100.02530678696402:(0.03520175882521048,0.1689626876062166)},xshift=0.03520175882521048cm,yshift=0.1689626876062166cm]  ({0.08319873700947768*(-1-(\x)^2)/(1-(\x)^2)},{0.1447370102073152*(-2)*(\x)/(1-(\x)^2)})  
--  (-1.93185,0.51764) circle (1.cm)
-- 
plot[samples=50,domain=-0.91:0.81,rotate around={-100.02530678696402:(0.03520175882521048,0.1689626876062166)},xshift=0.03520175882521048cm,yshift=0.1689626876062166cm]  ({0.08319873700947768*(1+(\x)^2)/(1-(\x)^2)},{0.1447370102073152*2*(\x)/(1-(\x)^2)}) --(1,0) circle (0.5cm)
;
\draw (-1.9318516525781364,-0.517638090205042) circle (1.cm);
\draw [samples=50,domain=-0.81:0.91,rotate around={-79.99135503997014:(0.0289278371406075,-0.17128698002763867)},xshift=0.0289278371406075cm,yshift=-0.17128698002763867cm] plot ({0.08317903326741531*(-1-(\x)^2)/(1-(\x)^2)},{0.14410135610422933*(-2)*(\x)/(1-(\x)^2)});
\draw [samples=50,domain=-0.91:0.81,rotate around={-79.99135503997014:(0.0289278371406075,-0.17128698002763867)},xshift=0.0289278371406075cm,yshift=-0.17128698002763867cm] plot ({0.08317903326741531*(1+(\x)^2)/(1-(\x)^2)},{0.14410135610422933*2*(\x)/(1-(\x)^2)});

\fill[fill opacity=0.3] (1,0) circle (0.5cm) --
plot[samples=50,domain=-0.81:0.91,rotate around={-79.99135503997014:(0.0289278371406075,-0.17128698002763867)},xshift=0.0289278371406075cm,yshift=-0.17128698002763867cm]  ({0.08317903326741531*(-1-(\x)^2)/(1-(\x)^2)},{0.14410135610422933*(-2)*(\x)/(1-(\x)^2)})
--
(-1.9318516525781364,-0.517638090205042) circle (1.cm)
-- 
plot [samples=50,domain=-0.91:0.81,rotate around={-79.99135503997014:(0.0289278371406075,-0.17128698002763867)},xshift=0.0289278371406075cm,yshift=-0.17128698002763867cm]  ({0.08317903326741531*(1+(\x)^2)/(1-(\x)^2)},{0.14410135610422933*2*(\x)/(1-(\x)^2)})--(1,0) circle (0.5cm)
;
\end{tikzpicture}\\
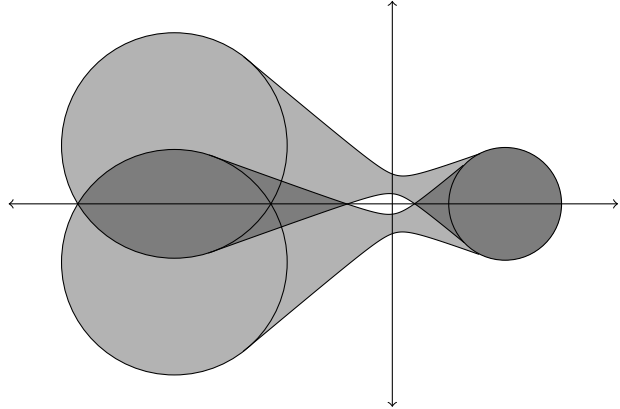
\textbf{Figure 10.} The product set $(K(1,2e^{i\frac{11\pi}{12}})\cup K(1,2e^{-i\frac{11\pi}{12}}))\cdot D(1,\frac{1}{2})$
is not simply connected.
\end{center}
\end{example}

\section{Additional results and further research}

We have to assume compactness in most of our results. One may wonder what happen
if we relax this assumption. The following example shows that
without the end points, the product of two line segments may not be
star-shaped.

\begin{example}  Let $K_1 = K_2$ be the line segment joining
$1+i$ and $1-i$ without the end points.
Then $K_1K_2$ has no star center.
\end{example}

\it
Verification. \rm  Note that the closure of $K_1K_2$ equals
$S = K(1+i,1-i)K(1+i,1-i)$ has a unique star-center 2. The set $K_1 K_2$
is obtained from $S$ by removing the line segments $K(2,2i)$ and $K(2,-2i)$.
The only point in the closure can reach all the points in $K_1K_2$ is 2, but it is not
in $K_1K_2$. So, $K_1K_2$ is not star-shaped.
\qed

Recall that an extreme point of a compact convex set $S \subseteq \IC$ is an element in
$S$ that cannot be written as the mid-point of two different elements in $S$.
If $S$ is a polygon (with interior) then its vertexes are the extreme points.
We can extend Theorem \ref{3.3} to the following.

\begin{theorem} \label{6.2}
Let $K_1, K_2 \subseteq \IC$ be compact convex sets.
Then $K_1K_2$ is star-shaped if and only if there is $p \in K_1K_2$
such that $K(p,ab) \subseteq K_1K_2$ for any extreme points $a \in K_1$ and
$b \in K_2$.
\end{theorem}

\it Proof. \rm If $K_1K_2$ is star-shaped, then a star center $p \in K_1K_2$ satisfies
$K(p,c) \subseteq K_1K_2$ for any $c \in K_1K_2$.
Now, suppose there is $p \in K_1K_2$ satisfying $K(p,ab)\subseteq K_1K_2$
for any extreme points $a \in K_1$ and $b \in K_2$.
Let $\mu = \mu_1\mu_2$ with $\mu_1 \in K_1, \mu_2 \in K_2$.
By the Caretheodory theorem $\mu_1 \in K(a_1, a_2, a_3)$ and
$\mu_2 \in K(b_1, b_2, b_3)$ for some extreme points
$a_1, a_2, a_3 \in K_1$ and $b_1, b_2, b_3 \in K_2$. (Some of the
$a_i$'s may be the same, and also some of the $b_i$'s may be the same.)
Suppose $p = p_1p_2$ with $p_1 \in K_1$ and $p_2 \in K_2$.
Then $p_1 \in K(a_4,a_5,a_6)$ and
$p_2 \in K(b_4, b_5,b_6)$
for some extreme points
$a_4, a_5, a_6 \in K_1$ and $b_4, b_5, b_6 \in K_2$.
By Theorem \ref{3.3},
$K(p,\mu_1\mu_2) \subseteq K(a_1, \dots, a_6)K(b_1, \dots, b_6) \subseteq K_1K_2$.
Thus, $p$ is a star center of $K_1K_2$ \qed

Another observation is the following extension of Proposition \ref{1.0}(b).
Note that we do not need to impose
compactness conditions on $K_1$ or $K_2$.

\begin{proposition}
Suppose $K_1 \subseteq \IC$ is star-shaped with 0 as a star center.
Then for any non-empty subset $K_2 \subseteq \IC$, the set $K_1K_2$
is star-shaped with 0 as a star center.
\end{proposition}

\it Proof. \rm Let $p = p_1 p_2 \in K_1 K_2$ with $p_1 \in K_1, p_2 \in K_2$.
Then $K(0,p) = K(0,p_1)\{p_2\} \subseteq K_1K_2$. \qed

There are other interesting questions deserve further research.
We mention a few of them in the following.

\begin{itemize}
\item[{\bf P1}] Find  necessary and sufficient conditions on $K_1$ and $K_2$
so that $K_1K_2$ is convex or star-shaped.
\end{itemize}

In the context of numerical range if $A \in M_2$, then $W(A)$ is an elliptical disk.
So, it is also of interest to study the following.

\begin{itemize}
\item[{\bf P2}]  Let $K_1, K_2$ be two elliptical disks. Determine conditions
on $K_1, K_2$ so that $K_1K_2$ is star-shaped or convex.
\end{itemize}

One may also consider the following.

\begin{itemize}

\item[{\bf P3}] Characterize those elliptical disks $K_1$ such that $K_1K_2$ is
star-shaped for all compact convex set $K_2$.
\end{itemize}

More generally, one may consider  the following.

\begin{itemize}
\item[{\bf P4}] Characterize those compact convex sets $K_1$ such that
$K_1K_2$ is convex or star-shaped for any compact convex set $K_2$.
\end{itemize}

In connection to Problem P4,  we have shown that if $K_1$ is a close line segment or
a close circular disk,  then $K_1K_2$ is star-shaped for any compact convex set $K_2$.
These results are are also connected to Problem P3 because
a line segment and a circular disk can be viewed as elliptical disks.

It is also interesting to study the Minkowski product of $s$ (convex) sets
$K_1, \dots, K_s$. The study will be more challenging. As pointed out in
\cite{Karol}, the set $K_1 \cdots K_s$ may not be simply connected in general.
Nevertheless, our results in Section 5 and Proposition 6.2 imply the
following.

\begin{proposition} Suppose $K_1, \dots, K_s \subseteq \IC$.

{\rm 1.} If any one of the sets $K_1, \dots, K_s$ is star-shaped with 0 as a star center, then
$K_1 \cdots K_s$ is star-shaped with 0 as a star center.

{\rm 2.} Suppose there is a nonzero number $\mu_1$ such that
$\mu_1 K_1$ is a circular disk center at 1 with radius $r < 1$.

{\rm (2.a)}
If there is $\mu \in \IC$ such that
$\mu K_2 \cdots K_s \subseteq \mu_1 K_1$ for $i = 2, \dots, s$, then
$K_1 \cdots K_r$ is star-shaped with $(\mu_1\mu)^{-1}(1-r^2)$ as a star center.

{\rm (2.b)}
If there is $\mu \in \IC$ such that
$\mu K_2 \cdots K_s \subseteq \{z \in \IC: \Re(z) \ge 1\}$ for $i = 2, \dots, s$, then
$K_1 \cdots K_r$ is star-shaped with $(\mu_1\mu)^{-1}$ as a star center.

\end{proposition}

It is also interesting to study the following problem.

\begin{itemize}
\item[{\bf P5}]
Characterize those compact (convex) sets $K$ such that
$K^2$ is convex or star-shaped.
\end{itemize}

\section*{Acknowledgment}

Li is an affiliate member of the Institute for Quantum Computing, University of Waterloo;
he is an honorary professor of the University of Hong Kong and the Shanghai University.
His research was supported by  USA NSF grant DMS 1331021  and 
Simons Foundation Grant  351047.
The research of Wang was supported by the Ministry of Science and Technology of the Republic
of China under project MOST 104-2918-I-009-001.

\end{document}